\documentclass[12pt]{article}

\usepackage{amssymb}
\usepackage{graphics} 
\usepackage{url}
\usepackage{hyperref}
\usepackage[authoryear,round]{natbib}
\usepackage{amssymb,color}
\usepackage{graphics,amsmath}

\begin{document}

\newenvironment {proof}{{\noindent\bf Proof.}}{\hfill $\Box$ \bigskip}

\newtheorem{theorem}{Theorem}[section]
\newtheorem{lemma}[theorem]{Lemma}
\newtheorem{condition}[theorem]{Condition}
\newtheorem{proposition}[theorem]{Proposition}
\newtheorem{remark}[theorem]{Remark}
\newtheorem{hypothesis}[theorem]{Hypothesis}
\newtheorem{corollary}[theorem]{Corollary}
\newtheorem{example}[theorem]{Example}
\newtheorem{definition}[theorem]{Definition}

\renewcommand {\theequation}{\arabic{section}.\arabic{equation}}

\def \non {{\nonumber}}
\def \noin {{\noindent}}
\def \hat {\widehat}
\def \tilde {\widetilde}
\def \E {\mathbb{E}}
\def \P {\mathbb{P}}
\def \R {\mathbb{R}}
\def \N {\mathbb{N}}
\def \1 {\mathbf{1}}

\title{\large {\bf Localization for constrained martingale problems and } 
\\{\bf optimal conditions for uniqueness of reflecting diffusions in 2-dimensional domains}}
                                                       
\author{\begin{tabular}{ll}                               
Cristina Costantini & Thomas G. Kurtz \\
Dipartimento di Economia & Department of Mathematics\\
and INdAM local unit & and Department of Statistics\\
Universit\`a di Chieti-Pescara & University of Wisconsin - Madison \\            
v.le Pindaro 42 & 480 Lincoln Drive  \\                                                    
65127 Pescara, Italy & Madison, WI  53706-1388, USA \\                    
c.costantini@unich.it & kurtz@math.wisc.edu    \\
& \url{http://www.math.wisc.edu/~kurtz/}  \\                   
\end{tabular}}

\date{April, 14 2023}

\maketitle

We prove existence and uniqueness for semimartingale reflecting diffusions in 2-dimensional piecewise smooth domains with varying, oblique directions of reflection on each "side", under geometric, easily verifiable conditions. Our conditions are optimal in the sense that, in the case of a convex polygon with constant direction of reflection on each side, they reduce to the conditions of Dai and Williams (1996), which are necessary for existence of Reflecting Brownian Motion. Moreover our conditions allow for cusps.

Our argument is based on a new localization result for constrained martingale problems which holds quite generally: as an additional example, we show that it holds for diffusions with jump boundary conditions. 
\vspace{.1in}


\vspace{.1in}


\newpage
\section{Introduction}\label{sectionintro}
Reflecting diffusions in nonsmooth domains have been 
studied since the early 1980s. 
Despite this long history, there is no 
general existence and uniqueness result in the literature 
for curved, piecewise smooth domains or cones, not even 
under the restriction that the process be a 
semimartingale, and not even in dimension 2. 
This notwithstanding the fact that there are significant 
applications, for instance in stochastic networks (see 
e.g. \cite{KKLW09} or \cite{KW12}). 

Exhaustive results exist only for normal reflection 
(\citet{Tan79}, \citet{Sai87}, \citet{BH91}, 
\citet{Bas96}, \citet{DBT93}, etc.), for Brownian motion in an orthant with constant direction of reflection on each face 
(\citet{HR81}, \citet{RW88}, \citet{TW93}, etc.), for Brownian motion in a 2-dimensional wedge 
with constant direction of reflection on each side 
(\citet{VW85}, \citet{Wil85}), 
for Brownian motion in a smooth cone with radially constant 
direction of reflection (\citet{KW91}) 
and for semimartingale reflecting Brownian motion 
in a convex polyhedral domain with constant direction of reflection on each 
face (\citet{DW96}). In the case of a simple polyhedral 
domain, the assumptions of \citet{DW96} are necessary 
for existence of a semimartingale Brownian motion (see also \citet{RW88} for the orthant case.)

For a piecewise smooth domain with 
varying, oblique direction of reflection on each ``face'', the best 
available result is \cite{DI93}. Unfortunately, the 
\cite{DI93} result is proved under a condition that is not 
easy to verify and leaves 
out many very natural examples. (See e.g. Remark 
\ref{re:DI}.) In fact, the \cite{DI93} condition does not 
reduce to the assumptions of \cite{DW96} in the case of a 
polyhedral domain. 

More recently, existence and uniqueness of a 
semimartingale reflecting diffusion has been proved by 
\cite{CK18} in a $2$-dimensional cusp with varying, oblique 
directions of reflection on each ``side'' and by \cite{CK22} 
in a $d$-dimensional domain with one singular point that 
near the singular point can be approximated by a smooth 
cone, with varying, oblique direction of reflection on the 
smooth part of the boundary.  In the cusp case, even 
starting at the cusp, with probability one, the process 
never hits it again.  In contrast, in the case when the 
domain can be approximated by a cone, the process can 
hit the singular point infinitely many times.  Therefore 
the study of this case requires a new ergodic theorem 
for inhomogeneous subprobability transition kernels.  The conditions 
under which the above results are proved are geometric in 
nature and easily verifiable. 
A quite general existence result for piecewise 
smooth domains in $\R^d$, even with cusp like points, has been obtained in \cite{CK19}, 
leaving the question of uniqueness.  

In dimension two, 
piecewise smooth domains look locally like smooth 
domains or like domains with one singular point.  
Consequently, by a localization argument, one should be 
able to exploit the results of \cite{CK18} and \cite{CK22} 
to give conditions for 
uniqueness of semimartingale reflecting diffusions. 
In this paper we carry out this program. 
The conditions we find (Conditions \ref{D} and \ref{G}; 
see also Remark \ref{re:C2}) 
are geometric and easy to verify 
and of course allow for cusps and for points where the 
boundary is smooth but the direction of reflection has a 
discontinuity.  
The same conditions allow to apply the results of \cite{CK19} to obtain 
existence as well. 
They are optimal in the sense that for a 
polygonal domain with constant direction of reflection on 
each side they reduce to the conditions of 
\cite{DW96} (Proposition \ref{th:DW} and Remark 
\ref{re:DW}.)

The existence proof in \cite{CK19} makes use of the 
equivalence between solutions of a stochastic differential 
equation with reflection (SDER) and natural solutions of 
the corresponding constrained martingale problem (CMP), 
proved in the same paper. 
CMPs were introduced in \cite{Kur90} and \cite{Kur91}) and further 
studied in \cite{KS01}, \cite{CK15} and \cite{CK19}. 
Here we exploit the equivalence between SDERs and 
CMPs also to localize the uniqueness problem for the SDER.  

In Section \ref{sectionlocliz}, we introduce 
CMPs stopped at the exit 
from an open set and show that, under a quite general 
condition, 
uniqueness holds for the natural solution of a CMP in a 
given domain if and only if it holds for the natural 
solution of the    
CMP stopped at the exit from each open set belonging to 
an open covering of the 
domain. This result holds for general CMPs in arbitrary 
dimension and 
is of independent interest. CMPs may be used to define 
not only reflecting diffusions, but also, for instance, diffusions with 
Wentzell boundary conditions and Markov processes with jump 
boundary conditions (see Section 7 of \cite{CK19}.) 
As an example of application of our localization result to 
other processes besides reflecting diffusions, we show 
that the condition we require is typically satisfied 
also by diffusions with jump boundary 
conditions (Remark \ref{re:ext}.) Since the proofs of the 
results of Section \ref{sectionlocliz} are somewhat 
technical, they are postponed to Appendix 
\ref{sectionloclizpf}. 

In Section \ref{sectionrefdif}, we combine the 
above localization results with the uniqueness results in 
\cite{CK18} and \cite{CK22} to obtain global uniqueness 
for the natural solution of the CMP corresponding to an 
SDER in a piecewise smooth domain in $\R^2$, with varying, 
oblique direction of reflection on each ``side''. As 
mentioned above, existence follows from \cite{CK19}. 
By the equivalence between natural solutions of the CMP and 
solutions of the SDER, existence and uniqueness transfer to the 
SDER. 
Although most of the work of this section consists 
in verifying the assumptions of \cite{CK19}, \cite{CK18} 
and \cite{CK22}, this verification is nontrivial. In 
particular, if the boundary has cusps, in order to apply 
the results of \cite{CK19} one needs to use the fact 
that the domain admits infinitely many representations 
and to construct a suitable representation. 

A more detailed discussion of the contents is provided 
at the beginning of each section. 

We will use the following notation. 
$\subseteq$ and $\supseteq$ will denote inclusion, while $\subset$ and $\supset$ will denote strict inclusion. 
For a finite set $F$, $|F|$ will denote the cardinality of $F$. 
For a metric space 
$E$, ${\cal B}(E)$ will denote the $\sigma$-algebra of Borel sets and $
{\cal P}(E)$ 
will denote the set of probability measures on $(E,{\cal B}(E))$; 
for $E_0\subseteq E$, $\overline {E_0}$ will denote the closure of $
E_0$. 
For a stochastic process $Z$, ${\cal F}^Z_t:=\sigma (Z(s),\,s\leq 
t)$ and 
${\cal F}^Z_{t^{+}}:=\cap_{s>t}\,{\cal F}^Z_s$;  
Finally the superscript $T$ denotes 
the transpose of a matrix and $B_r(0)$ denotes a ball in $\R^d$of 
radius $r$ and center the origin .

\section{Localization for constrained martingale problems}\label{sectionlocliz}

\setcounter{equation}{0}

Let $E$ be a compact metric space, $E_0$ be an open subset of 
$E$, and let $A\subseteq C(E)\times C(E)$ with $(1,0)\in A$.  
Let ${\cal U}$ also be a compact metric space, let $\Xi$ be a 
closed subset of 
$(E-E_0)\times {\cal U}$ and assume that, for every $x\in E-E_0$, 
there is some $u\in {\cal U}$ such that 
$(x,u)\in\Xi$. Let $B\subseteq C(E)\times C(\Xi )$ 
with $(1,0)\in B$, ${\cal D}:={\cal D}(A)\cap {\cal D}(B)$ and assume $
{\cal D}$ is dense 
in $C(E)$. The intuition is that $A$ is the generator for a 
process in $E$ and that $B$ determines controls that 
constrain the process to remain in $E_0$ or, more precisely, 
in $\overline {E_0}$.

Let ${\cal L}_{{\cal U}}$ be the space of Borel measures $\mu$ on $
[0,\infty )\times {\cal U}$ 
such that $\mu ([0,t]\times {\cal U})<\infty$ for all $t>0$.  ${\cal L}_{
{\cal U}}$ is
topologized so that  
$\mu_n\in {\cal L}_{{\cal U}}\rightarrow\mu\in {\cal L}_{{\cal U}}$ if and only if 
\[\int_{[0,\infty )\times {\cal U}}f(s,u)\mu_n(ds\times du)\rightarrow
\int_{[0,\infty )\times {\cal U}}f(s,u)\mu (ds\times du)\]
for all continuous $f$ with compact support in $[0,\infty )\times 
{\cal U}$. It is 
possible to define a metric on ${\cal L}_{\cal U}$ that induces the above 
topology and makes ${\cal L}_{\cal U}$ into a complete, separable metric 
space. Also let ${\cal L}_{\Xi}$ be defined analogously. For any 
${\cal L}_{\cal U}$-valued (${\cal L}_{\Xi}$-valued) random variable $L$, for each $
t\geq 0$, 
$L([0,t]\times\cdot )$ is a random measure on ${\cal U}$ ($\Xi$). We will 
occasionally use the notation $L(t):=L([0,t]\times\cdot )$. 

For a nondecreasing path $l_0\in D_{[0,\infty )}[0,\infty )$ with $
l_0(0)=0$, 
we define 
\begin{equation}(l_0)^{-1}(t):=\inf\{s\geq 0:\,l_0(s)>t\},\label{eq:tinv}\end{equation}
where we adopt the usual convention that the infimum of the empty set 
is $\infty$. Of course, if $l_0$ is strictly increasing $(l_0)^{-
1}$ is just 
the inverse of $l_0$. In addition, for every path $y\in D_E[0,\infty 
)$ 
or $y\in D_{[0,\infty )}[0,\infty )$ 
such that $\lim_{t\rightarrow\infty}y(t)$ exists, we will use the notation 
$y(\infty ):=\lim_{t\rightarrow\infty}y(t)$. 

The controlled martingale problem for $(A,E_0,B,\Xi )$, 
the constrained martingale problem for $(A,E_0,B,\Xi )$
and natural solutions of the constrained 
martingale problem for $(A,E_0,B,\Xi )$ have been introduced and studied in 
\cite{Kur90}, \cite{Kur91}, \cite{CK15} and \cite{CK19}. 
Here, given an open subset $U$ of $E$, we introduce the notions of {\em stopped controlled }
{\em martingale problem\/} for $(A,E_0,B,\Xi ;U)$, 
and of {\em natural solution of the }
{\em stopped constrained martingale problem\/} for $(A,E_0,B,\Xi 
;U)$ 
and study their relations with the corresponding 
unstopped objects. 
Our main goals are COrollary 
\ref{th:fromunst} and Theorem \ref{th:fromst}, which correspond 
to Theorems 4.6.1 and 4.6.3 of \cite{EK86} for martingale problems. 
A natural solution of the stopped constrained martingale problem 
for $(A,E_0,B,\Xi ;U)$ is obtained by time-changing a 
solution of the stopped controlled martingale problem 
for $(A,E_0,B,\Xi ;U)$ (see below for precise definitions): Roughly speaking, 
in order to transfer the results of Section 4.6 
of \cite{EK86} to constrained martingale 
problems, what we need is to be able to exchange 
the ''stopping'' and the ''time-changing''. 

Note that the set $E$ here corresponds to 
$\overline {E_0}\cup F_1$ in \cite{CK19} and that for Lemma \ref{th:paste} 
below we do not need Condition 3.5 c) of \cite{CK19}.

\begin{definition}\label{def:sclmp}
Let $Y^U$ be a process in $D_E[0,\infty )$, 
$\lambda_0^U$ be a nonnegative, nondecreasing process such that 
\begin{equation}\lambda_0^U(t)=\int_{[0,t]}{\bf 1}_{\overline {E_
0}}(Y^U(s))d\lambda_0^U(s)\quad a.s.,\label{eq:lam0supp}\end{equation}
and $\Lambda_1^U$ be a ${\cal L}_{{\cal U}}$-valued random variable such that  
\begin{equation}\lambda_1^U(t):=\Lambda_1^U([0,t]\times {\cal U})
=\int_{[0,t]\times {\cal U}}{\bf 1}_{\Xi}(Y^U(s),u)\Lambda_1^U(ds
\times du).\label{eq:Lam1supp}\end{equation}
Define
\begin{equation}\theta^U:=\inf\{t\geq 0:\,Y^U(t)\notin U\mbox{\rm \ or }
Y^U(t-)\notin U\}.\label{eq:extscl}\end{equation}
$(Y^U,\lambda_0^U,\Lambda_1^U)$ is a solution of the  {\em stopped, controlled }
{\em martingale problem\/} for $(A,E_0,B,\Xi ;U)$ if
\[(Y^U,\lambda_0^U,\Lambda_1^U)(t)=(Y^U,\lambda_0^U,\Lambda_1^U)(
t\wedge\theta^U),\quad\forall t\geq 0\quad a.s.,\]
\[\lambda_0^U(t)+\lambda_1^U(t)=t\wedge\theta^U,\quad\forall t\geq 
0\quad a.s.,\]
and 
\begin{equation}f(Y^U(t))-f(Y^U(0))-\int_0^tAf(Y^U(s))d\lambda_0^
U(s)-\int_{[0,t]\times {\cal U}}Bf(Y^U(s),u)\Lambda_1^U(ds\times 
du)\label{eq:sclmp}\end{equation}
is a $\{{\cal F}_t^{Y^U,\lambda_0^U,\Lambda_1^U}\}$-martingale for all $
f\in {\cal D}$. Since 
$(\ref{eq:sclmp})$ is right continuous, it is also a 
$\{{\cal F}_{t^{+}}^{Y^U,\lambda_0^U,\Lambda_1^U}\}$-martingale. 

For $U=E$, $(Y^U,\lambda_0^U,\Lambda_1^U)=(Y,\lambda_0,\Lambda_1)$ is a solution of the {\em controlled }
{\em martingale problem\/} for $(A,E_0,B,\Xi )$. 
\end{definition}

\begin{remark}\label{re:tost}
Note that, in general, $Y^U(t)$, in particular $Y^U(0)$, may 
take values outside $\overline U$. 
 
Let $(Y,\lambda_0,\Lambda_1)$ be a solution 
of the controlled martingale problem for $(A,E_0,B,\Xi )$. 
Then, setting 
\begin{equation}\theta :=\inf\{t\geq 0:\,Y(t)\notin U\mbox{\rm \ or }
Y(t-)\notin U\},\label{eq:extcl}\end{equation}
$(Y,\lambda_0,\Lambda_1)(\cdot\wedge\theta )$ is a solution of the 
stopped controlled martingale problem 
for $(A,E_0,B,\Xi ;U)$. 
\end{remark}
\vskip.2in

\begin{theorem}\label{th:paste}
Suppose that for every $\nu\in {\cal P}(E)$ there exists a solution 
of the controlled martingale problem for $(A,E_0,B,\Xi )$ with 
initial distribution $\nu$. 

Then, for every solution 
of the stopped controlled martingale problem 
for $(A,E_0,B,\Xi ;U)$, $(Y^U,\lambda_0^U,\Lambda_1^U)$, there exists a 
solution $(Y,\lambda_0,\Lambda_1)$ 
of the controlled martingale problem for $(A,E_0,B,\Xi )$ such 
that, with $\theta$ defined by $(\ref{eq:extcl})$, 
$\big(Y^U,\lambda_0^U,\Lambda_1^U,\theta^U\big)$ has the same distribution as 
$\big(Y(\cdot\wedge\theta ),\lambda_0(\cdot\wedge\theta ),\Lambda_1(\cdot
\wedge\theta ),\theta \big)$. 
\end{theorem}

\begin{proof}
See Appendix \ref{sectionloclizpf}. 
\end{proof}

\begin{definition}\label{def:scmp}
A process $X^U$ in $D_{\overline {E_0}}[0,\infty )$ is a solution of the {\em stopped }
{\em constrained martingale problem\/} for $(A,E_0,B,\Xi ;U)$ if there 
exists a ${\cal L}_{\Xi}$-valued random variable $\Lambda^U$ such that, setting 
\begin{equation}\tau^U:=\inf\{t\geq 0:\,X^U(t)\notin U\mbox{\rm \ or }
X^U(t-)\notin U\},\label{eq:extsc}\end{equation}
$(X^U,\Lambda^U)$ satisfies
\[(X^U,\Lambda^U)(t)=(X^U,\Lambda^U)(t\wedge\tau^U)\quad a.s.\]
and 
\begin{equation}f(X^U(t))-f(X^U(0))-\int_0^{t\wedge\tau^U}Af(X^U(
s))ds-\int_{[0,t]\times\Xi}Bf(x,u)\Lambda^U(ds\times dx\times du)\label{eq:scmp}\end{equation}
is a $\{{\cal F}_t^{X^U,\Lambda^U}\}$-local martingale for all $f
\in {\cal D}$. 
Since $(\ref{eq:scmp})$ is right continuous, it is also a 
$\{{\cal F}_{t^{+}}^{X^U,\Lambda^U}\}$-local martingale. 

For $U=E$, $X^U=X$ is a solution of the 
{\em constrained martingale problem\/} for $(A,E_0,B,\Xi )$ and we 
write $\Lambda^U=\Lambda$. 
\end{definition}

\begin{definition}\label{def:snat}
A solution $X^U$ of the stopped constrained martingale problem 
for $(A,E_0,B,\Xi ;U)$ is {\em natural}, if there exists a 
solution $(Y^U,\lambda_0^U,\Lambda_1^U)$ of the stopped controlled martingale 
problem, 
with the property that the event $\{\theta^U=\infty ,\;\lim_{s\rightarrow
\infty}\lambda_0^U(s)<\infty \}$ has 
zero probability, such that 
\[X^U(t)=Y^U((\lambda_0^U)^{-1}(t))\]
and 
\begin{equation}\Lambda^U([0,t]\times C):=\int_{[0,(\lambda_0^U)^{
-1}(t)]\times {\cal U}}{\bf 1}_C(Y^U(s),u)\Lambda_1^U(ds\times du
),\quad C\in {\cal B}(\Xi ),\quad a.s..\label{eq:Lam}\end{equation}
(Note that, a.s., if $\lim_{s\rightarrow\infty}\lambda_0^U(s)=t_0
<\infty$, then $\theta^U<\infty$ and 
$(\lambda_0^U)^{-1}(t)=\infty$ for all $t\geq t_0$, so that, for $
t\geq t_0$, 
$Y^U((\lambda_0^U)^{-1}(t))=Y^U(\infty )=Y^U(\theta^U)$.)

A solution $X$ of the constrained martingale problem 
for $(A,E_0,B,\Xi )$ is {\em natural}, if there exists a 
solution $(Y,\lambda_0,\Lambda_1)$ of the controlled martingale 
problem such that 
\[X(t)=Y((\lambda_0)^{-1}(t))\]
and 
\[\Lambda ([0,t]\times C):=\int_{[0,(\lambda_0)^{-1}
(t)]\times {\cal U}}{\bf 1}_C(Y(s),u)\Lambda_1(ds\times du),\quad 
C\in {\cal B}(\Xi ),\quad a.s..\]
\end{definition}

\begin{definition}\label{def:snat!}
Uniqueness holds for natural solutions of 
the stopped constrained martingale problem 
for $(A,E_0,B,\Xi ;U)$ 
(the constrained martingale problem for $(A,E_0,B,\Xi ;U)$)
if any two solutions with the same 
initial distributions have the same 
distribution on $D_{\overline {E_0}}[0,\infty )$. 
\end{definition}

In the sequel we assume the following condition on 
the controlled martingale problem for $(A,E_0,B,\Xi )$ and 
the open set $U$. 

\begin{condition}\label{hp:ext}\hfill 
\begin{itemize}
\item[(i)]For each $\nu\in {\cal P}(E)$ there exists a 
solution $(Y,\lambda_0,\Lambda_1)$ of the controlled martingale problem for $
(A,E_0,B,\Xi )$ 
with initial distribution $\nu$.
\end{itemize}

For each solution $(Y,\lambda_0,\Lambda_1)$ of the controlled martingale 
problem for $(A,E_0,B,\Xi )$:
\begin{itemize}
\item[(ii)]
\[\lim_{t\rightarrow\infty}\lambda_0(t)=\infty\quad a.s..\]
\item[(iii)]There exists a sequence of 
$\bigg\{{\cal F}^{Y,\lambda_0,\Lambda_1}_{\lambda_0^{-1}(t)}\bigg
\}$- stopping times $\{\gamma_n\}$ such that $\gamma_n\rightarrow
\infty$ a.s. and 
$\E[\lambda_0^{-1}(\gamma_n)]<\infty$ for each $n$. 

\item[(iv)]
For $X(t):=Y(\lambda_0^{-1}(t))$, $\tau$ defined as 
\begin{equation}\tau :=\inf\{t\geq 0:\,X(t)\notin U\mbox{\rm \ or }
X(t^{-})\notin U\}\quad a.s.\label{eq:extc}\end{equation}
and $\theta$ defined by $(\ref{eq:extcl})$,   
\[\lambda_0^{-1}(\tau )=\theta\quad a.s..\]
\end{itemize}
\end{condition}
\vskip.2in

\begin{remark}\label{re:CMP}
(i) and (ii) of Condition \ref{hp:ext} are a) and b) of Condition 3.5 
of \cite{CK19}. Together with (iii), they ensure that $X$, 
defined as in (iv), is a natural solution of the constrained 
martingale problem of $(A,E_0,B,\Xi )$: See 
Theorem 3.6 of \cite{CK19}. 
\end{remark}

\begin{proposition}\label{th:inc}
Suppose Condition \ref{hp:ext} (i) is verified. 
If each solution of the controlled martingale 
problem for $(A,E_0,B,\Xi )$ satisfies $\lambda_0(t)>0$ for all $
t>0$ 
a.s., then $\lambda_0$ is strictly increasing a.s. for each solution, 
and Condition \ref{hp:ext} is verified for every open set 
$U$.  
\end{proposition}

\begin{proof}
See Appendix \ref{sectionloclizpf}
\end{proof}

\begin{remark}\label{re:ext}
The controlled martingale problems corresponding to 
reflecting diffusions will usually satisfy the assumptions 
of Proposition \ref{th:inc} (e.g. see Lemma 6.8 of 
\cite{CK19}). However there are significant 
examples of controlled martingale problems 
for which Condition \ref{hp:ext} is verified for a large 
class of open sets $U$ 
although the assumptions of Proposition \ref{th:inc} are 
not satisfied. For instance, this is the case for 
diffusions with jump 
boundary conditions. Let $E_0$ be a bounded domain in $\R^d$ 
with smooth boundary, $E$ be a compact set in $\R^d$ such that 
$\overline {E_0}\subseteq\stackrel {\circ}E$, where $\overline {E_
0}$ and $\stackrel {\circ}E$ denote the closure of $E_0$ and the interior of $
E$ in the 
topology of $\R^d$ respectively. Consider the operator 
\[Af(x):=b(x)\cdot\nabla f(x)+\frac 12\mbox{\rm tr}((\sigma\sigma^
T)(x)D^2f(x)),\]
where $\sigma\sigma^T$ is uniformly positive definite on $\overline {
E_0}$, $b$ and 
$\sigma$ are continuous and vanish outside of an open 
neighborhood of $\overline {E_0}$ whose closure is included in $\stackrel {
\circ}E$. 
Let ${\cal U}:=\{1\}$, $\Xi :=(E-E_0)\times {\cal U}$ and $B$ be defined by 
\[Bf(x,1)=Bf(x):=\int (f(y)-f(x))p(x,dy),\]
where $p$ is a transition function on $E$, $p(x,\cdot )$ is 
continuous as a function from $E$ into ${\cal P}(E)$ and, for all 
$x\in E$, 
\[p(x,E_0)=1.\]
Then the controlled martingale problem for $(A,E_0,B,\Xi )$ 
satisfies (i), (ii) and (iii) of Condition \ref{hp:ext}: 
see Section 7.1 of \cite{CK19}, and note that, 
under the above assumptions, Lemma 3.1 of 
\cite{CK19} applies, so that (iii) holds with $\gamma_n:=n$. 

Intuitively if $(Y,\lambda_0,\Lambda_1)$ is a solution of the controlled martingale 
problem for $(A,E_0,B,\Xi )$, $Y$ behaves like a 
diffusion with generator $A$ till it reaches $\partial E_0$; it 
stays at the exit point for a unit exponential time and 
then it jumps into $E_0$ and starts behaving like a diffusion 
again. The corresponding natural solution of the 
constrained martingale problem for $(A,E_0,B,\Xi )$ defined in 
(iv) of Condition \ref{hp:ext} 
behaves in the same way except that it jumps 
instantaneously. In particular both $Y$ and $X$ stay in $\overline {
E_0}$ 
for all times and $Y$ ($X$) jumps at a time $t$ if and only if 
$Y(t^{-})\in\partial E_0$ ($X(t^{-})\in\partial E_0$). 

If $Y(0)\in\partial D$, $Y$ will stay at $Y(0)$ 
for a unit exponential time $\rho$ and $\lambda_0(t)=0$ for 
$0<t\leq\rho$, therefore the assumption of Proposition 
\ref{th:inc} is not satisfied. 
However, let $U$ be an open set of $\R^d$ with smooth boundary, 
such that $\overline U\subseteq\stackrel {\circ}E$ and that, denoting by $\mbox{\rm Leb}$ the surface Lebesgue measure 
on $\partial U$, $\mbox{\rm Leb}(\partial U\cap\partial E_0)=0$. Then, with $
\theta$ and $\tau$ as in (iv) of Condition 
\ref{hp:ext}, the probability that $Y(\theta^{-})$ belongs to $\partial 
U\cap\partial E_0$ 
is zero. It follows that, almost surely, either $Y(\theta^{-})\in
\partial U\cap E_0$, so 
that $Y(\theta )=Y(\theta^{-})\in E_0$, or $Y(\theta^{-})\in U$ and $
Y(\theta )\notin U$, so that 
$Y(\theta^{-})\in\partial E_0$ and $Y(\theta )\in E_0$. In both cases 
$\lambda_0$ is strictly increasing in a right neighborhood ot $\theta$, so 
that $\lambda_0^{-1}(\lambda_0(\theta ))=\theta$. Moreover $Y(\theta 
)\notin U$ implies $\tau =\lambda_0(\theta )$, so 
that (iv) of Condition \ref{hp:ext} is satisfied. 

Processes of this type have been considered in a variety 
of settings, for example \cite{DN90,SS94}.  Semigroups 
corresponding to processes with nonlocal boundary 
conditions of this type have been considered in 
\cite{AKK16}.  Related models are considered in 
\cite{MR85}.

\end{remark}
\vskip.2in

\begin{theorem}\label{th:snatext}
Under Condition \ref{hp:ext}, 
for every natural solution $X^U$ of the stopped constrained 
martingale problem for $(A,E_0,B,\Xi ;U)$, 
there exists a natural solution $X$ of the constrained martingale problem 
for $(A,E_0,B,\Xi )$ such that, with $\tau$ defined by $(\ref{eq:extc})$, 
$X(\cdot\wedge\tau )$ has the same distribution as $X^U(\cdot )$. 
\end{theorem}

\begin{proof}
See Appendix \ref{sectionloclizpf}
\end{proof}

\begin{corollary}\label{th:fromunst}
Under Condition \ref{hp:ext}, if uniqueness holds for natural 
solutions of the constrained martingale problem for 
$(A,E_0,B,\Xi )$, then it holds for natural solutions of 
the stopped constrained martingale problem for 
$(A,E_0,B,\Xi ;U)$. 
\end{corollary}

\begin{proof}
The assertion follows immediately from Theorem  
\ref{th:snatext}. 
\end{proof}

\begin{theorem}\label{th:fromst} 
Suppose there exist open subsets 
 $U_k\subseteq E$, $k=1,2,...$, with $E=\bigcup_{k=1}^{\infty}U_k$, such that, for 
each $k$, $(A,E_0,B,\Xi )$ and $U_k$ satisfy Condition \ref{hp:ext}
and uniqueness holds for natural solutions of 
the stopped, constrained martingale problem for 
$(A,E_0,B,\Xi ;U_k)$. 
Then uniqueness holds for natural solutions of the 
constrained martingale problem for $(A,E_0,B,\Xi )$.  
\end{theorem}

\begin{proof}
See Appendix \ref{sectionloclizpf}. 
\end{proof}

\section{Existence and uniqueness of reflecting diffusions in a 2-dimensional, piecewise smooth domain}\label{sectionrefdif}

\setcounter{equation}{0}

In this section, first we formulate our assumptions on the 
domain where the reflecting diffusion takes values and 
on the directions of reflection and 
compare them with the assumptions of the most general previous results, 
namely the results of \cite{DI93} (Remark \ref{re:DI}) and 
\cite{DW96} (Proposition \ref{th:DW}). In particular, in 
the case of a convex polygon with constant direction of 
reflection on each side,  
our assumptions are equivalent to those of 
\cite{DW96}, which are necessary for 
existence of a reflecting Brownian motion: in this sense 
our assumptions are optimal (Remark \ref{re:DW}). 

Next we prove that the two definitions of a 
semimartingale reflecting 
diffusion as a solution of a stochastic differential 
equation with reflection and as a natural solution of a 
constrained martingale problem are equivalent (Theorem 
\ref{th:equiv}) and prove existence of a reflecting 
diffusion (Theorem \ref{th:exist}). Both these results 
follow immediately from the results of Section 6 of 
\cite{CK19} once one has verified that the assumptions of 
Section 6 of \cite{CK19} are satisfied (Lemma 
\ref{th:sec6.19}: however, in particular at a cusp point, 
this verification is nontrivial and requires to construct a suitable 
representation of the domain. 

Finally, we show that uniqueness holds for the constrained 
martingale problem stopped at the exit from a 
neighborhood of each corner, both when the corner is a 
cusp (Lemma \ref{th:cuspunq}) and when it is not (Lemma 
\ref{th:coneunq}): this amounts essentially to verifying that the 
assumptions of \cite{CK18} and \cite{CK22}, respectively, 
are satisfied, but, again, this is nontrivial. Corollary 
\ref{th:fromunst} is also needed here. 
Uniqueness for the global constrained martingale problem 
then follows immediately from Theorem \ref{th:fromst} 
and transfers to the corresponding stochastic differential 
equation by Theorem \ref{th:equiv}.

We consider a domain $D$ satisfying the following 
condition. 

\begin{condition}\label{D}\hfill 
\begin{itemize}
\item[(i)]$D$ is a bounded domain that 
admits the representation 
\begin{equation}D=\bigcap_{i=1}^mD^i,\label{eq:repr}\end{equation}
where, for $i=1,...,m$, $D^i$ is a bounded domain defined as 
\[D^i:=\{x:\psi^i(x)>0\},\quad\psi^i\in {\cal C}^1(\R^2),\quad\inf_{}{}_{
x:\,\psi^i(x)=0}\,|\nabla\psi^i(x)|>0,\]
and 
\[\overline D=\bigcap_{i=1}^m\overline {D^i}.\]
The representation is minimal in the 
sense that, for $j=1,...,m$, 
\[D\subset\bigcap_{i\neq j}D^i,\]
where $\subset$ denotes strict inclusion. 

For $x\in\partial D^i$, we denote by $n^i(x)$ the unit, inward normal 
to $D^i$ at $x$, i.e. $n^i(x):=\frac {\nabla\psi^i(x)}{|\nabla\psi^
i(x)|}$.
\item[(ii)] For $x^0\in\bigcup_{i=1}^m\partial D^i$ and
\begin{equation}I(x^0):=\{i:\,x^0\in\partial D^i\},\label{eq:indexset}\end{equation}
$ $the set $\{x\in\bigcup_{i=1}^m\partial D^i:\,|I(x^0)|>1\}$ is finite. We call a 
point $x^0\in\partial D$ such that $|I(x^0)|>1$ a {\em corner\/} and  
assume $|I(x^0)|=2$ at every corner. 
 
\item[(iii)]Let $x^0$ be a corner and $I(x^0)=\{i,j\}$. 

If $n^i(x^0)\neq -n^j(x^0)$
  (then we say that $x^0$ is a {\em cone point\/}), 
\[\limsup_{x\in\partial D^l-\{x^0\},\,x\rightarrow x^0}\frac {|n^
l(x)-n^l(x^0)|}{|x-x^0|}<\infty ,\quad\limsup_{x\in\partial D^l-\{
x^0\},\,x\rightarrow x^0}\frac {|n^l(x^0)\cdot (x-x^0)|}{|x-x^0|^
2}<\infty ,\]
for $l=i,j$. 

If $n^j(x^0)=-n^i(x^0)$ (then we say that $x^0$ is a {\em cusp point\/}), 
$D\cap B_r(x^0)$ is connected for all $r>0$ small enough, and  
\[\lim_{x\in\partial D^i\cap\partial D-\{x^0\},\,z\in\partial D^j
\cap\partial D-\{x^0\},\,|(x-z)\cdot n^i(x^0)|=|x-z|,\,x,z\rightarrow 
x^0}\frac {(x-x^0)\cdot n^i(x^0)}{(x-z)\cdot n^i(x^0)}=L,\]
for some finite $L$.
\end{itemize}
\end{condition}
\vskip.2in

\begin{remark}\label{re:repr}
A piecewise ${\cal C}^1$ domain $D$ admits infinitely many 
representations $(\ref{eq:repr})$, 
and it may be that some representations verify all 
assumptions in Condition \ref{D} and others do not. 
In all our results we only need that there exists a 
representation that verifies Condition \ref{D}. 
It may be convenient to use more than one 
representation with different properties (see Lemma 
\ref{th:sec6.19}). 
\end{remark}
\vskip.2in

Define the inward normal cone at $x^0\in\partial D$ as 
\begin{equation}N(x^0):=\bigg\{n:\,\liminf_{x\in\overline D-\{x^0
\},\,x\rightarrow x^0}\frac {(x-x^0)}{|x-x^0|}\cdot n\geq 0\bigg\}
.\label{eq:Ncone}\end{equation}
For $I(x^0)=\{i,j\}$, if $x^0$ is a cone point, clearly 
$N(x^0)$ is the closed, convex cone generated by $n^i(x^0)$ and 
$n^j(x^0)$. If $x^0$ is a cusp point, by the assumption that 
$D\cap\partial B_r(0)$ is connected for all $r>0$ small enough, there exists one and only one 
unit vector $\tau (x^0)$ such that 
\begin{equation}\tau (x^0)\cdot n^i(x^0)=0\quad\mbox{\rm and}\quad\lim_{
x\in\overline D-\{x^0\},\,x\rightarrow x^0}\frac {\tau (x^0)\cdot 
(x-x^0)}{|x-x^0|}=1.\label{eq:tau}\end{equation}
Then 
\begin{equation}N(x^0)=\{u\in\R^2:\,u\cdot\tau (x^0)\geq 0\}.\label{eq:cuNcone}\end{equation}
\vskip.2in

\begin{remark}\label{re:C2}
Let $x^0$ be a corner, $I(x^0)=\{i,j\}$, and suppose 
$\psi^i,\psi^j\in\!{\cal C}^2(\R^2)$. Then, if $x^0$ is a cone point, Condition 
\ref{D} (iii) is always verified; if $x^0$ is a cusp point 
Condition \ref{D} (iii) is verified if 
\[\tau (x^0)\cdot\bigg(\frac {D^2\psi^j(x^0)}{|\nabla\psi^j(x^0)|}
+\frac {D^2\psi^i(x^0)}{|\nabla\psi^i(x^0)|}\bigg)\tau (x^0)\neq 
0.\]
\end{remark}
\vskip.2in

The set of possible directions of reflection on the 
boundary of $D$ is defined by vector fields $g^i:\R^2\rightarrow\R^
2$, 
$i=1,...,m$, $g^i$ of unit length on $\partial D^i$. For $x^0\in\partial 
D$, define 
\begin{equation}G(x^0):=\bigg\{\sum_{i\in I(x^0)}\eta_ig^i(x^0),\,
\eta_i\geq 0\bigg\}.\label{eq:Gcone}\end{equation}

\begin{condition}\label{G}\hfill
\begin{itemize}
\item[(i)]For $i=1,...,m$, $g^i$ is a Lischitz continuous vector 
field such that 
\[\inf_{x\in\partial D^i}g^i(x)\cdot n^i(x)>0.\]
\item[(ii)]For every $x^0\in\partial D$, there exists a unit vector 
$e(x^0)\in N(x^0)$ such that 
\[e(x^0)\cdot g>0,\quad\forall g\in G(x^0)-\{0\}.\]
\end{itemize}
\end{condition}

\begin{remark}\label{re:DI}
As mentioned in the Introduction, the best 
result available in the literature for a piecewise smooth domain
with varying directions of reflection on each ''face'' is \cite{DI93}. 
A very simple example that shows how the \cite{DI93} 
assumptions may not be satisfied is the following. 
Let $D^1$ be the unit ball centered at $(1,0)$, and let $D$ be its 
intersection with the upper half plane. Of course $D$ 
can be represented as $D:=D^1\cap D^2$, where 
$D^2$ is a bounded ${\cal C}^1$ domain. 
Let $n^i$, $i=1,2$, denote the unit, inward normal to $D^i$, and 
\[g^i(x)\equiv\left[\begin{array}{cc}
\cos(\vartheta )&\sin(\vartheta )\\
-\sin(\vartheta )&\cos(\vartheta )\end{array}
\right]n^i(x),\qquad\vartheta\mbox{\rm \ a constant angle, }\frac {
\pi}4\leq\vartheta <\frac {\pi}2.\]
Then, at $x^0=(0,0)$ and at $x^0=(2,0)$, it can be proved by contradiction that 
there is no convex compact set that satisfies $(3.7)$ of 
\cite{DI93}. Conditions \ref{D} and \ref{G} are instead 
satisfied. 
\end{remark}
\vskip.2in

In the case when $D$ is a convex 
polygon and the direction of reflection is constant on 
each side, Condition \ref{G} coincides with 
the assumptions of \cite{DW96}. This is an immediate 
consequence of the following lemma, which rephrases the 
assumptions of \cite{DW96}. The lemma holds in general for convex 
polyhedrons in $\R^d$. 

Let 
\begin{equation}D:=\bigcap_{i=1}^m\{x\in\R^d:\,x\cdot n^i>b_i\},\quad i=1,...,m
,\label{eq:DWdom}\end{equation}
where $n^1,...,n^m$ are distinct unit vectors, $b_1,...,b_m$ 
are real numbers, and the above representation is 
minimal, that is, for each $j=1,...,m$,
\begin{equation}D\subset\bigcap_{i\neq j}\{x\in\R^d:\,x\cdot n^i>
b_i\},\label{eq:DWmin}\end{equation}
where $\subset$ denotes strict inclusion. 

Assumption 1.1 of \cite{DW96} is formulated in terms of 
{\em maximal\/} subsets of the set of indeces $\{1,...,m\}$, defined as 
follows: ${\bf K}\subseteq \{1,...,m\}$ is maximal if and only if 
${\bf K}\neq\emptyset$, $F_{{\bf K}}:=\{x\in\overline D:\,x\cdot 
n^i=b_i,\,\forall i\in {\bf K}\}\neq\emptyset$ and, for every 
${\bf K}'\supset {\bf K}$, $F_{{\bf K}'}\subset F_{{\bf K}}$ (where $
\supset$ ad $\subset$ $ $denote strict inclusion). 

\begin{lemma}\label{th:DWmax}
${\bf K}\subseteq \{1,...,m\}$ is maximal if and only if ${\bf K}
=I(x^0)$ 
for some $x^0\in\partial D$. 
\end{lemma}

\begin{proof}
For ${\bf K}=\{1,...,m\}$, being maximal is equivalent to 
$F_{{\bf K}}\neq\emptyset$, that is ${\bf K}=I(x^0)$ for some $x^
0\in\partial D$. 

For ${\bf K}\subset \{1,...,m\}${\bf ,} ${\bf K}$ is 
maximal if and only if for every $j\notin {\bf K}$ there exists 
$x^j\in\overline D$ such that $x^j\cdot n^i=b_i$ for all $i\in {\bf K}$, $
x^j\cdot n^j>b_j$. 
Then the fact that ${\bf K}=I(x^0)$ is maximal for every 
$x^0\in\partial D$ is immediate. To see that the converse holds, 
let ${\bf K}$ be maximal and set  
\[x^0:=\frac 1{m-|{\bf K}|}\sum_{j\in \{1,...,m\}-{\bf K}}x^j.\]
Then $x^0\in\overline D$ and 
\[x^0\cdot n^i=b_i,\quad\forall i\in {\bf K},\]
\[x^0\cdot n^i=\frac 1{m-|{\bf K}|}\big(x^i\cdot n^i+\sum_{j\in \{
1,...,m\}-{\bf K},\,j\neq i}x^j\cdot n^i\big)>b_i,\quad\forall i\in 
\{1,...,m\}-{\bf K},\]
that is ${\bf K}=I(x^0)$. 
\end{proof}

\begin{proposition}\label{th:DW}
Let $D\subseteq\R^2$ be defined by $(\ref{eq:DWdom})$ and be 
bounded, and let 
$g^i$, $i=1,...,m$, be constant unit vectors. 

Then $D$ satisfies Condition \ref{D}. $D$ and $g^i$, $i=1,...,m$, satisfy 
Condition \ref{G} if and only if they satisfy Assumption 1.1 of 
\cite{DW96}. 
\end{proposition}

\begin{proof}
Verifying that $D$ satisfies Condition \ref{D} is immediate. 
In particular, in this case the minimality assumption 
$(\ref{eq:DWmin})$ implies that $1\leq |I(x^0)|\leq 2$ for every 
$x^0\in\partial D$. 

In dimension 2 every polyhedron 
is simple (see Definition 1.4 of \cite{DW96}), therefore, 
by Proposition 1.1 of \cite{DW96}, Assumption 1.1 of 
\cite{DW96} reduces to assuming that, for each maximal 
${\bf K}$, there is a nonnegative linear combination $e:=\sum_{i\in 
{\bf K}}\eta_in^i$ 
such that $e\cdot g^{_j}>0$ for all $j\in {\bf K}$ (actually \cite{DW96} 
requires a positive linear combination, but of 
course the two requirements are equivalent). 
Since, by Lemma \ref{th:DWmax}, ${\bf K}$ is maximal if and only if $
{\bf K}=I(x^0)$ for some 
$x^0\in\partial D$, this is indeed Condition \ref{G} (ii). 
As the directions of reflection $g^i$ are constant, Condition \ref{G} (i) follows from (ii). 
\end{proof}

\begin{remark}\label{re:disc-dir}
Conditions \ref{D} and \ref{G} allow for boundary points 
$x^0$ at which the boundary is actually smooth, but 
the direction of reflection has a discontinuity, i.e. 
\[n^i(x^0)=n^j(x^0),\quad g^i(x^0)\neq g^j(x^0),\qquad i,j\in I(x^
0).\]
\end{remark}
\vskip.2in

Finally, we assume that the drift $b$ and the dispersion coefficient $
\sigma$ satisfy the following 
condition. 

\begin{condition}\label{A}\hfill 
\begin{itemize}
\item[(i)] $b:\R^2\rightarrow\R^2$ and $\sigma :\R^2\rightarrow\R^{
2\times 2}$ are Lipschitz continuous. 

\item[(ii)]For every corner $x^0$, $\sigma (x^0)$ is non singular. 
\end{itemize}
\end{condition}

In most of the literature, 
a semimartingale reflecting diffusion is defined as a solution of a 
{\em stochastic differential equation with reflection}. We 
recall the definition below, for the convenience of the 
reader. 

\begin{definition}\label{def:SDER}
Let $D$ be a bounded domain and, for $x\in\partial D$, let $G(x)$ be
 a closed, convex cone 
such that $\{(x,u)\in\partial D\times\partial B_1(0):\,u\in G(x)\}$ is closed.  Let 
$b:\R^2\rightarrow\R^2$ and $\sigma :\R^2\rightarrow\R^{2\times 2}$ be bounded, measurable 
functions, and $\nu\in {\cal P}(\overline D)$. 
A stochastic process $X$ is a solution of the stochastic 
differential equation with reflection in $\overline D$ with 
coefficients $b$ and $\sigma$, 
cone of directions of reflection $G$, and initial distribution 
$\nu$, if $X(0)$ has distribution $\nu$, there exist a 
standard Brownian motion $W$, a continuous, non 
decreasing process $\lambda$, and a process 
$\gamma$ with measurable paths, all defined on the same probability space as $
X$, 
such that $W(t+\cdot )-W(t)$ is independent of 
${\cal F}_t^{X,W,\lambda ,\gamma}$, for all $t\geq 0$, and the equation 
\begin{eqnarray}
&&X(t)=X(0)+\int_0^tb(X(s))ds+\int_0^t\sigma (X(s))dW(s)+\int_0^t
\gamma (s)\,d\lambda (s),\quad t\geq 0,\non\\
&&\quad\gamma (t)\in G(X(t)),\quad |\gamma (t)|=1,\quad d\lambda 
-a.e.,\quad t\geq 0,\label{eq:SDER}\\
&&X(t)\in\overline D,\quad\lambda (t)=\int_0^t{\bf 1}_{\partial D}
(X(s))d\lambda (s),\quad t\geq 0,\non\end{eqnarray}
is satisfied a.s..

Given an initial distribution $\nu\in {\cal P}(\overline D)$, {\em weak }
{\em uniqueness\/} or {\em uniqueness in distribution\/} holds if all 
solutions of $(\ref{eq:SDER})$ with
$P\{X(0)\in\cdot \}=\nu$ have the same distribution on 
$C_{\bar {D}}[0,\infty )$. 

A stochastic process $\tilde {X}$  
is a weak solution of $(\ref{eq:SDER})$ if there is a 
solution $X$ of (\ref{eq:SDER}) such that $\tilde {X}$ and $X$ have the 
same distribution.
\end{definition}\

\begin{remark}\label{re:DW}
When $D$ is a bounded, convex polyhedron in $\R^2$, and the direction of 
reflection is constant on each side, Propositions 1.1 and 1.2 of \cite{DW96} prove that if 
there exists a semimartingale reflecting Brownian 
motion (i.e. a weak solution of $(\ref{eq:SDER})$ with $b$ and 
$\sigma$ constant), then Assumption 1.1 of \cite{DW96} must be 
verified. On the other hand we have proved in 
Proposition \ref{th:DW} that, when specialized to this 
case, Condition \ref{G} coincides with Assumption 1.1 of 
\cite{DW96}. In this sense Condition \ref{G} is optimal. 
\end{remark}

In the following we exploit repeatedly the equivalence 
between the stochastic differential equation 
$(\ref{eq:SDER})$ and the constrained martingale 
problem for $(A,D,B,\Xi )$, where the state space is $E:=\overline 
D$, 
$A$ denotes the operator 
\begin{equation}{\cal D}(A):={\cal C}^2(\overline D),\quad Af(x):
=b(x)\cdot\nabla f(x)+\frac 12\mbox{\rm tr}((\sigma\sigma^T)(x)D^
2f(x)),\label{eq:A}\end{equation}
and 
\begin{eqnarray}&&{\cal U}:=\partial B_1(0),\quad\Xi :=\{(x,u)\in\partial 
D\times {\cal U}:\,u\in G(x)\}, \label{eq:B}\\
&&B:{\cal C}^2(\overline D)\rightarrow 
{\cal C}(\Xi ),\,\,\,Bf(x,u):=\nabla f(x)\cdot u.\nonumber
\end{eqnarray}

This equivalence is proved in general dimension $d$ 
in Section 6 of \cite{CK19} (Theorem 6.12), under quite 
general assumptions. In the next lemma we show that, 
under Conditions \ref{D}, \ref{G} and \ref{A}, the 
assumptions of Section 6 of \cite{CK19} are satisfied, or 
more precisely, that the domain $D$ admits a representation 
such that the assumptions of Section 6 of \cite{CK19} are 
verified (see Remark \ref{re:repr}.)

\begin{lemma}\label{th:sec6.19}
Assume Conditions \ref{D}, \ref{G} and \ref{A}. Then the 
domain $D$ admits a representation 
\[D=\bigcap_{i=1}^{\tilde {m}}\tilde {D}^i,\]
such that the assumptions of Section 6 of \cite{CK19} are 
verified. 
\end{lemma}

\begin{proof}
First of all note that the assumption of Section 6 of \cite{CK19} 
that the domains are simply connected is 
redundant: it is enough to assume that the domains are 
connected, as we are doing here. 

Let $x^0$ be a corner. We suppose, without loss 
of generality, that $x^0=0$, $I(0)=\{1,2\}$, and we 
write $n^1$ for $n^1(0)$ and $n^2$ for $n^2(0)$. 

If $0$ is a cone point, the normal cone $N(0)$ can be 
written in the form $(6.3)$ and Conditions 6.2 a) and b) of \cite{CK19} 
are verified. Condition \ref{G} (ii) implies that the matrix 
\[\left[\begin{array}{cc}
n^1\cdot g^1&n^2\cdot g^1\\
n^1\cdot g^2&n^2\cdot g^2\end{array}
\right]\]
is a completely-S matrix. Then its transpose 
is also completely-S (Lemma 3 of \cite{RW88}), so that, 
in particular, there exists $g\in G(0)$, $g=c_1g^1+c_2g^2$, 
$c_1,c_2>0$, such that $n^1\cdot g>0$, $n^2\cdot g>0$. Therefore, for each 
$n=\eta_1n^1+\eta_2n^2$, $\eta_1,\eta_2\geq 0$, $\eta_1+\eta_2>0$, $
n\cdot g>0$, i.e.  
$c_1n\cdot g^1+c_2n\cdot g^2>0$, which implies that $n\cdot g^1>0$ or 
$n\cdot g^2>0$, that is Condition 6.2 (c) of \cite{CK19} for 
$I=\{1,2\}$. Since Condition 6.2 (c) of \cite{CK19} is clearly 
satisfied for $I=\{1\}$ and $I=\{2\}$, it is verified for every 
$I\subseteq I(0)$.  

Now let $0$ be a cusp point and let $\tau =\tau (0)$ be the 
vector defined in $(\ref{eq:tau})$. Without loss of generality we 
can take $(\tau ,n^1)$ as the basis of the coordinate system.  
Let $r_0>0$ be small enough that $\overline {B_{r_0}(0)}$ contains 
no other corners than $0$. Then $D$ can be represented as 
\begin{equation}D=\Delta\cap\tilde {D}^1\cap\tilde {D}^2\cap\bigcap_{i\geq 3}D^
i,\label{eq:repr1}\end{equation}
($\bigcap_{i\geq 3}D^i=\R^2$ if $m=2$), with $\Delta$ a bounded domain 
with ${\cal C}^1$ boundary, such that 
\[\Delta\cap\overline {B_{r_0}(0)}=\{x\in\overline {B_{r_0}(0)}:\,
x\cdot\tau >0\},\quad\Delta\supseteq\overline D-\{0\},\]
and 
\[\tilde {D}^i:=\{x:\tilde{\psi}^i(x)>0\},\qquad i=1,2,\]
\[\tilde{\psi}^i(x_1,x_2):=\psi^i(|x_1|,x_2)\big[1-\chi\big(\frac 
2{r_0}(|x|-\frac {r_0}2)\big)\big]+\psi^i(x_1,x_2)\chi\big(\frac 
2{r_0}(|x|-\frac {r_0}2)\big),\qquad i=1,2,\]
where $\psi^i$ is the function defining $D^i$ and $\chi$ is a smooth, 
nondecreasing function such that 
$\chi (t)=0$ for $t\leq 0$, $\chi (t)=1$ for $t\geq 1$. 

Intuitively, we add an extra domain $\Delta$ and replace the 
function $\psi^i$, $i=1,2$, with a function $\tilde{\psi}^i$ that agrees with 
$\psi^i$ for $x_1\geq 0,$ but is symmetric 
with respect to $x_1$ in a neighborhood of $0$ 
With the addition of the extra domain $\Delta$, the normal cone 
$N(0)$ can be written in the form $(6.3)$ of \cite{CK19}. By defining the 
direction of reflection on $\partial\Delta$, $\gamma$, to be the inward normal 
direction, we have $\gamma (0)=\tau$, so that the cone of 
directions of reflection at $0$, $G(0)$, does not change. 
However, by the symmetry of the functions $\tilde{\psi}^i$, $i=1,
2$, 
now $\tilde {{\cal I}}(0)$, defined by $(6.9)$ of \cite{CK19} for the 
representation $(\ref{eq:repr1})$, is  
\[\tilde {{\cal I}}(0)=\big\{\{0\},\{1\},\{2\},\{0,1\},\{0,2\}\big
\},\]
and Condition 6.2 (c) of \cite{CK19} is satisfied at $0$.  

By iterating the above construction for each cusp point of 
$\partial D$, we obtain a representation of $D$ that satisfies the 
assumptions of Section 6 of \cite{CK19}.   
\end{proof}

\begin{theorem}\label{th:equiv}
Every solution of $(\ref{eq:SDER})$ 
is a natural solution of the constrained martingale problem 
for $(A,D,B,\Xi )$ defined by $(\ref{eq:A})$-$(\ref{eq:B})$.

Conversely every natural solution 
of the constrained martingale problem for $(A,D,B,\Xi )$ 
is a weak solution of $(\ref{eq:SDER})$.
\end{theorem}

\begin{proof}\ By Lemma \ref{th:sec6.19}, this is just 
a special case of Theorem 6.12 of 
\cite{CK19}. Note that a solution of $(\ref{eq:SDER})$ as 
defined in Definition \ref{def:SDER} is 
called a weak solution in \cite{CK19}.  
\end{proof}

\begin{theorem}\label{th:exist}
Under Conditions \ref{D}, \ref{G} and \ref{A}, 
for every initial distribution $\nu\in {\cal P}(\overline D)$, there exists a 
strong Markov solution of $(\ref{eq:SDER})$ 
with initial distribution $\nu$. 
\end{theorem}

\begin{proof}
By Lemma \ref{th:sec6.19}, this is just 
a special case of Theorem 6.13 of \cite{CK19}. 
\end{proof}

\begin{remark}\label{re:numapp}
Note that the construction of the solution of 
$(\ref{eq:SDER})$ provided in Section 6 of \cite{CK19} 
(Theorem 6.7 of \cite{CK19} and Lemma 1.1 of 
\cite{Kur90}) yields also a numerical approximation of 
the solution. 
\end{remark}

\begin{lemma}\label{th:coneunq}
Let $x^0\in\partial D$ be a cone point, $r_0$ be small enough that 
$\partial D\cap\overline {B_{r_0}(x^0)}$ contains no other corners and 
$U:=\overline D\cap B_{r_0}(x^0)$. Let $A$, $\Xi$ and $B$ be defined by 
$(\ref{eq:A})$-$(\ref{eq:B})$. 

Then, under Conditions \ref{D}, \ref{G} and \ref{A}, uniqueness holds for 
natural solutions of the stopped constrained martingale 
problem for $(A,D,B,\Xi ;U)$. 
\end{lemma}

\begin{proof}
We suppose, without loss of generality, that $x^0=0$, 
$I(0)=\{1,2\}$, and we write $n^1$, $g^1$, $n^2$, $g^2$ for $n^1(
0)$, $g^1(0)$, 
etc.. 

Let $\tilde {D}\subseteq D$ be a bounded domain with boundary of class 
${\cal C}^1$ at every point except $0$, such that 
$\overline {\tilde {D}}\cap\overline {B_{r_0}(0)}=\overline D\cap\overline {
B_{r_0}(0)}$ and denote by $\tilde {n}(x)$ the unit, 
inward normal to $\tilde {D}$ at $x\in\partial\tilde {D}-\{0\}$. 
Let $\tilde {G}(x):=\{\eta\tilde {g}(x),\,\eta\geq 0\}$ for $x\in
\partial\tilde {D}-\{0\}$, where 
$\tilde {g}^{}:\R^2-\{0\}\rightarrow\R^2$ is some locally Lipschitz continuous vector 
field, of unit length on $\partial\tilde {D}-\{0\}$, 
such that $\tilde {g}(x)\cdot\tilde {n}(x)>0$ for $x\in\partial\tilde {
D}-\{0\}$ and 
$\tilde {G}(x)=G(x)$ for $x\in\big(\partial\tilde {D}\big)\cap\overline {
B_{r_0}(0)}-\{0\}$. Set $ $
\[\tilde {G}(0):=G(0),\qquad\tilde{\Xi }:=\{(x,u)\in\partial\tilde {
D}\times {\cal U}:\,u\in\tilde {G}(x)\}.\]
Let 
\[\tilde {{\cal K}}:=\{u\in\R^2:\,u\cdot n^1>0,\,u\cdot n^2>0\},\]
if $n^1\neq n^2$, and  
\[\tilde {{\cal K}}:=\{u\in\R^2:\,u\cdot n^1>0\},\]
if $n^1=n^2$. Then it can be checked by elementary 
computations that Condition \ref{D} implies that 
$\tilde {D}$ and $\tilde {{\cal K}}$ satisfy Conditions 3.1 (i) and (ii) of \cite{CK22}. 
Conditions 3.3 (i), (ii) and (iv) of \cite{CK22} also follow immediately from 
Condition \ref{G} (i) and (ii). 

As for Condition 3.3 (iii) of \cite{CK22}, if $n^1=n^2$ then 
$N(0)=\{\eta n^1,\,\eta\geq 0\}$ and 
Condition \ref{G} (ii) says that $G(0)-\{0\}\subseteq\tilde {{\cal K}}$. 
If $n^1\neq n^2$, by the argument already used in the proof of 
Lemma \ref{th:sec6.19}, Condition \ref{G} (ii) implies that the matrix 
\[\left[\begin{array}{cc}
n^1\cdot g^1&n^1\cdot g^2\\
n^2\cdot g^1&n^2\cdot g^2\end{array}
\right]\]
is a completely-S matrix, which in particular implies 
that there is $g\in G(0)$ such that $n^1\cdot g>0$. $n^2\cdot g>0$, i.e. 
$g\in G(0)\cap\tilde {{\cal K}}$. 

Therefore, by Theorem 3.25 of \cite{CK22}, uniqueness holds for 
natural solutions of the constrained martingale 
problem for $(A,\tilde {D},B,\tilde{\Xi })$. Moreover, it is shown in the 
proof of Theorem 3.23 of \cite{CK22} that, for each 
$\nu\in {\cal P}(\overline {\tilde {D}})$, there exists a solution of the controlled 
martingale problem for $(A,\tilde {D},B,\tilde{\Xi })$ with initial 
distribution $\nu$. Together with Lemma 3.16 of 
\cite{CK22} and 
Proposition \ref{th:inc}, this ensures that Condition \ref{hp:ext} 
is verified by $(A,\tilde {D},B,\tilde{\Xi })$ and $U$. 

Then Corollary \ref{th:fromunst} yields that uniqueness 
holds for natural solutions of the stopped constrained 
martingale problem for $(A,\tilde {D},B,\tilde{\Xi };U)$. A solution $
X^U$ of the 
the stopped constrained martingale problem for 
$(A,D,B,\Xi ;U)$ is not necessarily a 
solution of the stopped constrained 
martingale problem for $(A,\tilde {D},B,\tilde{\Xi };U)$ because its initial 
distribution might charge $\overline D\cap \big( \overline {\tilde {D}}\big)^c$. However if $X^U$ and 
$\tilde {X}^U$ are two solutions of the 
the stopped constrained martingale problem for 
$(A,D,B,\Xi ;U)$ with the same initial distribution, 
\[Z^U(t):=\left\{\begin{array}{ll}
X^U(t),\,t\geq 0,&\mbox{\rm if }X^U(0)\in U,\\
z^0,\,t\geq 0,&\mbox{\rm if }X^U(0)\notin U,\end{array}
\right.\qquad\tilde {Z}^U(t):=\left\{\begin{array}{ll}
\tilde {X}^U(t),\,t\geq 0,&\mbox{\rm if }\tilde {X}^U(0)\in U,\\
z^0,\,t\geq 0,&\mbox{\rm if }\tilde {X}^U(0)\notin U,\end{array}
\right.\]
where $z^0$ is some fixed point in $\overline {\tilde {D}}-U$, are two solutions 
of the stopped constrained 
martingale problem for $(A,\tilde {D},B,\tilde{\Xi };U)$ with the same 
initial distribution. Therefore $Z^U$ and $\tilde {Z}^U$ have the same distribution and so do 
$X^U$ and $\tilde {X}^U$.
\end{proof}

\begin{lemma}\label{th:cuspunq}
Let $x^0\in\partial D$ be a cusp point, $r_0$ be small enough that 
$\partial D\cap\overline {B_{r_0}(x^0)}$ contains no other corners and 
$U:=\overline D\cap B_{r_0}(x^0)$. Let $A$, $\Xi$ and $B$ be defined by 
$(\ref{eq:A})$-$(\ref{eq:B})$.

Then, under Conditions \ref{D}, \ref{G} and \ref{A}, uniqueness holds for 
natural solutions of the stopped constrained martingale 
problem for $(A,D,B,\Xi ;U)$.
\end{lemma}

\begin{proof}
Suppose, without loss of generality, that $x^0=0$, 
$I(0)=\{1,2\}$. We will write $n^1$, $g^1$, $n^2$, $g^2$ for $n^1
(0)$, $g^1(0)$, 
etc.. 

Let $\tilde {D}$, $\tilde {G}(x)$, $x\in\partial\tilde {D}-\{0\}$, and $
\tilde{\Xi}$ be as in the proof 
of Lemma \ref{th:coneunq}, in particular $\tilde {G}(0):=G(0).$
Let $\tau =\tau (0)$ be the vector in $(\ref{eq:tau})$ and take $
(\tau ,n^1)$ 
as the basis of the coordinate system. By the implicit 
function theorem there exist $r_1>0$, $r_2>0$, $r_1^2+r_2^2\leq r_
0^2$, and continuosly 
differentiable functions $\varphi^1$ and $\varphi^2$ defined on $
[-r_1,r_1]$, 
with values in $[-r_2,r_2]$, such that $\varphi^1(0)=\varphi^2(0)
=0$ 
and, for $(x_1,x_2)\in [-r_1,r_1]\times [-r_2,r_2]$,  
\[\psi^1(x_1,x_2)>0\quad\Leftrightarrow\quad x_2>\varphi^1(x_1),\qquad
\psi^1(x_1,x_2)=0\quad\Leftrightarrow\quad x_2=\varphi^1(x_1),\]
\[\psi^2(x_1,x_2)>0\quad\Leftrightarrow\quad x_2<\varphi^2(x_1),\qquad
\psi^2(x_1,x_2)=0\quad\Leftrightarrow\quad x_2=\varphi^2(x_1).\]
Then $\varphi^1$ and $\varphi^2$ satisfy Condition 2.1 of \cite{CK18}. 
In addition, taking into account $(\ref{eq:cuNcone})$, 
Condition \ref{G} ensures that $\tilde {g}$ satisfies Condition 2.3 of \cite{CK18}. 
Therefore Theorems 3.1, 4.1 and 4.7 of \cite{CK18}), together 
with Theorem \ref{th:equiv}, give 
uniqueness for natural solutions of the constrained martingale 
problem for $(A,\tilde {D},B,\tilde{\Xi })$. Moreover, in the proof of 
Theorem 4.1 of \cite{CK18} a solution of the controlled 
martingale problem for $(A,\tilde {D},B,\tilde{\Xi })$ with initial distribution 
the Dirac measure at $0$ is constructed and it is shown 
that, for that solution, $\lambda_0$ (denoted as $K_0$ there) is 
strictly increasing. Exactly the same arguments allow to 
construct a solution of the controlled 
martingale problem for $(A,\tilde {D},B,\tilde{\Xi })$ with an arbitrary initial 
distribution $\nu\in {\cal P}\big( \overline {\tilde {D}}\big)$ and to show that $
\lambda_0$ is 
strictly increasing for each solution of the controlled 
martingale problem for $(A,\tilde {D},B,\tilde{\Xi })$. Hence, by Proposition 
\ref{th:inc}, Condition \ref{hp:ext} is satisfied by 
$(A,\tilde {D},B,\tilde{\Xi })$ and $U$ and we can conclude as in the proof of 
Lemma \ref{th:coneunq}.
\end{proof}

\begin{theorem}\label{th:uniq}
Under Conditions \ref{D}, \ref{G} and \ref{A}, 
for every initial distribution $\nu\in {\cal P}(\overline D)$, uniqueness in 
distribution holds 
for solutions of $(\ref{eq:SDER})$ with initial distribution 
$\nu$. 
\end{theorem}

\begin{proof}
Let $A$, $\Xi$ and $B$ be defined by $(\ref{eq:A})$-$(\ref{eq:B})$.
By Lemma \ref{th:sec6.19}, 
$D$, $G$, $b$ and $\sigma$ satisfy the assumptions of Section 6 of \cite{CK19}, 
therefore Theorems 6.7 and Lemma 6.8 of \cite{CK19},  
together with Proposition \ref{th:inc}, ensure that  
Condition \ref{hp:ext} is satisfied by 
$(A,D,B,\Xi )$ and any open set $U$. 

Let $x^1$, $x^2$, ..., $x^M$ be the corners of $\overline D$, $r_
0>0$ 
be such that $x^h\notin\overline {B_{r_0}(x^k)}$ for $h\neq k$. Let 
$U^k:=\overline D\cap B_{r_0}(x^k)$, $k=1,...,M$, 
$U^{M+1}:=\overline D\cap\bigg(\bigcup_{k=1}^M\overline {B_{r_0/2}
(x^k)}\bigg)^c$. 

By Lemmas \ref{th:coneunq} and \ref{th:cuspunq}, uniqueness holds for 
natural solutions of the stopped constrained martingale 
problems for $(A,D,B,\Xi ;U^k)$, for $k=1,...,M$. As for the 
stopped constrained martingale 
problem for $(A,D,B,\Xi ;U^{M+1})$, one can consider a domain 
$\tilde {U}^{M+1}\subseteq D$ with ${\cal C}^1$ boundary, such that 
$\overline {\tilde {U}^{M+1}}\cap\bigg(\bigcup_{k=1}^MB_{r_0/2}(x^
k)\bigg)^c=\overline D\cap\bigg(\bigcup_{k=1}^MB_{r_0/2}(x^k)\bigg
)^c$ and a 
Lipschitz continuous direction of reflection $\tilde {g}^{M+1}$ on 
$\partial\tilde {U}^{M+1}$ such 
that $\tilde {G}^{M+1}(x):=\{\eta\tilde {g}^{M+1}(x),\,\eta\geq 0
\}=G(x)$ for 
$x\in\partial\tilde {U}^{M+1}\cap\bigg(\bigcup_{k=1}^MB_{r_0/2}(x^
k)\bigg)^c$, and argue as in 
Lemmas \ref{th:coneunq} and \ref{th:cuspunq}, but using 
Corollary 5.2 (Case 2) of 
\cite{DI93} and Theorem 6.12 of \cite{CK19}, to obtain that 
uniqueness holds for 
natural solutions of the stopped constrained martingale 
problems for $(A,D,B,\Xi ;U^{M+1})$. 

Then the assertion follows by Theorems \ref{th:fromst} 
and \ref{th:equiv}. 
\end{proof}

\appendix 
\section{Proofs of Section \ref{sectionlocliz}}\label{sectionloclizpf}
\renewcommand {\theequation}{A.\arabic{equation}} 
\setcounter{equation}{0} 
\vskip.2in

\noindent{\bf Proof of Theorem \ref{th:paste}}

The proof is a suitable modification of the proof of Lemma 4.5.16 of 
\cite{EK86}: 

Let $P^U$ denote the distribution of $(Y^U,\lambda_0^
U,\Lambda_1^U)$, 
$\nu$ denote the distribution of $Y^U(\theta^U)$ 
and $P$ denote the distribution of a solution 
of the controlled martingale problem for $(A,E_0,B,\Xi )$ with 
initial distribution $\nu$. 
let $Q$ be the probability measure on 
$D_E[0,\infty )\times C_{[0,\infty )}[0,\infty )\times {\cal L}_{\cal U}
\times [0,\infty ]\times D_E[0,\infty )\times C_{[0,\infty )}[0,\infty 
)\times {\cal L}_{\cal U}$ 
defined by 
\begin{equation}Q(D_1\times D_2):=\int_E\E^{P^U}\big[{\bf 1}_{D_1}
(\eta^1,l_0^1,L_1^1,\vartheta )\big|\eta^1(\vartheta )=y\big]\E^P\big
[{\bf 1}_{D_2}(\eta^2,l_0^2,L_1^2)\big|\eta^2(0)=y\big]\nu (dy)\label{eq:paste}\end{equation}
where $(\eta^1,l_0^1,L_1^1,\vartheta ,\eta^2,l_0^2,L_1^2)$ is the coordinate random 
variable in 
$D_E[0,\infty )\times C_{[0,\infty )}[0,\infty )\times {\cal L}_{\cal U}
\times [0,\infty ]\times D_E[0,\infty )\times C_{[0,\infty )}[0,\infty 
)\times {\cal L}_{\cal U}$, 
$D_1$ is a Borel subset of $D_E[0,\infty )\times C_{[0,\infty )}[
0,\infty )\times {\cal L}_{\cal U}\times [0,\infty ]$ 
and $D_2$ is a Borel subset of $D_E[0,\infty )\times C_{[0,\infty 
)}[0,\infty )\times {\cal L}_{\cal U}$. 

Define, for $t\geq 0$, $C\in {\cal B}({\cal U})$, 
\[Y(t):=\left\{\begin{array}{ll}
\eta^1(t),&t<\vartheta\\
\eta^2(t-\vartheta ),&t\geq\vartheta ,\end{array}
\right.\]
\[\lambda_0(t):=\left\{\begin{array}{ll}
l_0^1(t),&t<\vartheta\\
l_0^2(t-\vartheta )+l_0^1(\vartheta ),&t\geq\vartheta ,\end{array}
\right.\]
\[\Lambda_1([0,t]\times C):=\left\{\begin{array}{ll}
L_1^1([0,t]\times C)&t<\vartheta ,\\
L_1^2([0,t-\vartheta ]\times C)+L_1^1([0,\vartheta ]\times C),&t\geq
\vartheta ,\end{array}
\right.\]
\[\theta :=\vartheta\]
Then the distribution of $(Y,\lambda_0,\Lambda_1)(\cdot\wedge\theta 
)$ under $Q$ is 
$P^U$. In particular $\theta$ as defined above agrees $Q$-a.s. with 
$\theta$ as defined in $(\ref{eq:extcl})$. 

Let us show that $(Y,\lambda_0,\Lambda_1)$ is a solution of the controlled 
martingale problem for $(A,E_0,B,\Xi )$. To this end we need 
to show that, for arbitrary $0=t_0<t_1<...t_n<t_{n+1}$, denoting 
\[R:=f(Y(t_{n+1}))-f(Y(t_n))-\int_{t_n}^{t_{n+1}}Af(Y(s))d\lambda_
0(s)-\int_{[t_n,t_{n+1}]\times {\cal U}}Bf(Y(s),u)\Lambda_1(ds\times du)
,\]
it holds, for arbitrary continuous functions $h_k$ and $H_k$ and 
$C_k\in {\cal B}({\cal U})$,  
\[\E\big[R\prod_{k=1}^nh_k(Y(t_k))\,H_k(\lambda_0(t_k)-\lambda_0(
t_{k-1}),\Lambda_1((t_{k-1},t_k]\times C_k))\big]=0.\]
Observing that 
\[\lambda_0(t_k)-\lambda_0(t_{k-1})=\lambda_0(t_k\vee\theta )-\lambda_
0(t_{k-1}\vee\theta )+\lambda_0(t_k\wedge\theta )-\lambda_0(t_{k-
1}\wedge\theta ),\]
\[\Lambda_1((t_{k-1},t_k]\times C_k)=\Lambda_1((t_{k-1}\vee\theta 
,t_k\vee\theta ]\times C_k)+\Lambda_1((t_{k-1}\wedge\theta ,t_k\wedge
\theta ]\times C_k)\]
we see that we can replace 
$H_k(\lambda_0(t_k)-\lambda_0(t_{k-1}),\Lambda_1((t_{k-1},t_k]\times 
C_k))$ by the product 
\begin{eqnarray*}H_k^{\vee}(\lambda_0(t_k\vee\theta )-\lambda_0(t_{k-1}\vee\theta 
),\Lambda_1((t_{k-1}\vee\theta ,t_k\vee\theta ]\times C_k))\\
\times H_k^{
\wedge}(\lambda_0(t_k\wedge\theta )-\lambda_0(t_{k-1}\wedge\theta 
),\Lambda_1((t_{k-1}\wedge\theta ,t_k\wedge\theta ]\times C_k)),\end{eqnarray*}
where $H_k^{\vee}$ and $H_k^{\wedge}$ are arbitrary continuous functions 
such that $H_k^{\vee}(0,0)=H_k^{\wedge}(0,0)=1$. 
Analogously we can split $R$ as 
\begin{eqnarray*}
&&\qquad\qquad\qquad\qquad\qquad\qquad\qquad\qquad\qquad R=R^{\vee}
+R^{\wedge},\\
&&R^{\vee}:=
f(Y(t_{n+1}\vee\theta ))-f(Y(t_n\vee\theta ))\\
&&\qquad\qquad\qquad\qquad -\int_{t_
n\vee\theta}^{t_{n+1}\vee\theta}Af(Y(s))d\lambda_0(s)-\int_{(t_n\vee
\theta ,t_{n+1}\vee\theta ]\times {\cal U}}Bf(Y(s),u)\Lambda_1(ds\times 
du),\\
&&R^{\wedge}:=
f(Y(t_{n+1}\wedge\theta ))-f(Y(t_n\wedge\theta ))\\
&& \qquad\qquad\qquad\qquad -\int_{
t_n\wedge\theta}^{t_{n+1}\wedge\theta}Af(Y(s))d\lambda_0(s)-\int_{
(t_n\wedge\theta ,t_{n+1}\wedge\theta ]\times {\cal U}}Bf(Y(s),u)\Lambda_
1(ds\times du),\end{eqnarray*}
so that we reduce to proving that 
\begin{eqnarray}
&&\E^Q\big[R^{\vee}\prod_{k=1}^nh_k(Y(t_k))\,H_k^{\vee}(\lambda_0
(t_k\vee\theta )-\lambda_0(t_{k-1}\vee\theta ),\Lambda_1((t_{k-1}
\vee\theta ,t_k\vee\theta ]\times C_k))\non\\
&&\qquad\qquad\qquad H_k^{\wedge}(\lambda_0(t_k\wedge\theta 
)-\lambda_0(t_{k-1}\wedge\theta ),\Lambda_1((t_{k-1}\wedge\theta 
,t_k\wedge\theta ]\times C_k))\big]=0,\label{eq:vee}\end{eqnarray}
\begin{eqnarray}
&&\E^Q\big[R^{\wedge}\prod_{k=1}^nh_k(Y(t_k))\,H_k^{\vee}(\lambda_
0(t_k\vee\theta )-\lambda_0(t_{k-1}\vee\theta ),\Lambda_1((t_{k-1}
\vee\theta ,t_k\vee\theta ]\times C_k))\non\\
&&\qquad\qquad\qquad H_k^{\wedge}(\lambda_0(t_k\wedge\theta 
)-\lambda_0(t_{k-1}\wedge\theta ),\Lambda_1((t_{k-1}\wedge\theta 
,t_k\wedge\theta ]\times C_k))\big]=0.\label{eq:wed}\end{eqnarray}
Noting that 
\[R^{\wedge}=R^{\wedge}\,{\bf 1}_{\theta >t_n}\]
and that 
\[{\bf 1}_{\theta >t_n}\prod_{k=1}^nH_k^{\vee}(\lambda_0(t_k\vee\theta 
)-\lambda_0(t_{k-1}\vee\theta ),\Lambda_1((t_{k-1}\vee\theta ,t_k
\vee\theta ]\times C_k))\,\\
={\bf 1}_{\theta >t_n},\]
we see, by computations analogous to those of Lemma 
4.5.16 of \cite{EK86}, 
that the left hand side of $(\ref{eq:wed})$ equals zero. 

In order to see that $(\ref{eq:vee})$ is verified, define 
\[\vartheta_m:=\frac {[m\vartheta ]}m,\]
\begin{eqnarray*}
&&R^{\vee}_m:=f(\eta^2(t_{n+1}\vee\vartheta_m-\vartheta_m))-f(\eta^
2(t_n\vee\vartheta_m-\vartheta_m))-\int_{t_n\vee\vartheta_m-\vartheta_
m}^{t_{n+1}\vee\vartheta_m-\vartheta_m}Af(\eta^2(s))dl_0^2(s)\\
&&\qquad\qquad\qquad\qquad\qquad\quad\qquad\qquad -\int_{
(t_n\vee\vartheta_m-\vartheta_m,t_{n+1}\vee\vartheta_m-\vartheta_
m]\times {\cal U}}Bf(\eta^2(s),u)L_1^2(ds\times du),\end{eqnarray*}

and consider
\begin{eqnarray}
&&R^{\vee}_m\prod_{t_k<\vartheta_m}h_k(\eta^1(t_k))\,H_k^{\vee}(l_
0^1(t_k\vee\vartheta_m)-l_0^1(t_{k-1}\vee\vartheta_m),L_1^1((t_{k
-1}\vee\vartheta_m,t_k\vee\vartheta_m]\times C_k))\non\\
&&\qquad\qquad\qquad H_k^{\wedge}(l_0^1(t_k\wedge\vartheta_m)-l_0^
1(t_{k-1}\wedge\vartheta_m),L_1^1((t_{k-1}\wedge\vartheta_m,t_k\wedge
\vartheta_m]\times C_k))\non\\
&&\quad\;\;\prod_{t_k\geq\vartheta_m}h_k(\eta^2(t_k-\vartheta_m))\label{eq:veem}\\
&&\qquad\qquad\qquad H_k^{\vee}(l_0^2(t_k\vee\vartheta_m-\vartheta_
m)-l_0^2(t_{k-1}\vee\vartheta_m-\vartheta_m),\non \\
&&\qquad\qquad\qquad\qquad L_1^2((t_{k-1}\vee\vartheta_
m-\vartheta_m,t_k\vee\vartheta_m-\vartheta_m]\times C_k))\non\\
&&\qquad\qquad\qquad H_k^{\wedge}(l_0^2(t_k\wedge\vartheta_m-\vartheta_
m)-l_0^2(t_{k-1}\wedge\vartheta_m-\vartheta_m),\non \\
&&\qquad\qquad\qquad\qquad L_1^2((t_{k-1}\wedge
\vartheta_m-\vartheta_m,t_k\wedge\vartheta_m-\vartheta_m]\times C_
k)),\non\end{eqnarray}

Noting that 
\[\prod_{t_k<\vartheta_m}H_k^{\vee}(l_0^1(t_k\vee\vartheta_m)-l_0^
1(t_{k-1}\vee\vartheta_m),L_1^1((t_{k-1}\vee\vartheta_m,t_k\vee\vartheta_
m]\times C_k))=1,\]
\[\prod_{t_{k-1}\geq\vartheta_m}H_k^{\wedge}(l_0^2(t_k\wedge\vartheta_
m-\vartheta_m)-l_0^2(t_{k-1}\wedge\vartheta_m-\vartheta_m),L_1^2(
(t_{k-1}\wedge\vartheta_m-\vartheta_m,t_k\wedge\vartheta_m-\vartheta_
m]\times C_k))=1,\]

and that 
\[R^{\vee}_m=R^{\vee}_m\,{\bf 1}_{\vartheta_m<t_{n+1}},\]
we find, by computations analogous to those of Lemma 
4.5.16 of \cite{EK86}, that the expectation of $(\ref{eq:veem})$ under $
Q$ 
equals zero. 
Since $(\ref{eq:veem})$ converges pointwise and boundedly 
to 
$R^{\vee}\prod_{k=1}^nh_k(Y(t_k))\,H_k^{\vee}(\lambda_0(t_k\vee\theta 
)-\lambda_0(t_{k-1}\vee\theta ),\Lambda_1((t_{k-1}\vee\theta ,t_k
\vee\theta ]\times C_k))\,\\
H_k^{\wedge}(\lambda_0(t_k\wedge\theta )-\lambda_0(t_{k-1}\wedge\theta 
),\Lambda_1((t_{k-1}\wedge\theta ,t_k\wedge\theta ]\times C_k))$, 
$(\ref{eq:vee})$ is verified. 
\hfill $\Box$ \bigskip
\vskip.3in

\noindent {\bf Proof of Proposition \ref{th:inc}}

(ii) and the fact that $\lambda_0$ is 
strictly increasing follow from Lemmas 3.3 
and 3.4 of \cite{CK19}. In turn, the fact that $\lambda_0$ is 
strictly increasing immediately implies  
(iv). As in the proof of Corollary 3.9 of \cite{CK19}, 
(iii) is verified by 
\[\gamma_n:=\lambda_0(n).\]
\hfill $\Box$ \bigskip
\vskip.2in

\noindent {\bf Proof of Theorem \ref{th:snatext}}

Let $X^U(\cdot )=Y^U((\lambda_0^U)^{-1}(\cdot ))$ for some 
solution $(Y^U,\lambda_0^U,\Lambda_1^U)$ of the stopped controlled martingale problem 
for $(A,E_0,B,\Xi ;U)$, and let $(Y,\lambda_0,\Lambda_1)$ be the solution of the 
controlled martingale problem for $(A,E_0,B,\Xi )$ constructed 
in Theorem \ref{th:paste}.   Let $\theta$ be defined by 
$(\ref{eq:extcl})$. By Remark \ref{re:CMP}, $X(\cdot ):=Y(\lambda_
0^{-1}(\cdot ))$ 
is a natural solution of the constrained martingale problem for 
$(A,E_0,B,\Xi )$. Then, by Condition \ref{hp:ext} (iv), 
\[X(t\wedge\tau )=Y(\lambda_0^{-1}(t)\wedge\theta )=Y(\lambda_0(\cdot
\wedge\theta )^{-1}(t)\wedge\theta ),\quad t\geq 0,\]
and the assertion follows from the fact that the distribution of 
$Y(\lambda_0(\cdot\wedge\theta )^{-1}(\cdot )\wedge\theta )$ is the distribution of $
Y^U((\lambda_0^U)^{-1}(\cdot ))$, i.e. of 
$X^U(\cdot )$. 
\hfill $\Box$ \bigskip
\vskip.3in

\begin{lemma}\label{th:extc}
For each solution $(Y^U,\lambda_0^U,\Lambda_1^U)$ of the stopped controlled martingale 
problem for $(A,E_0,B,\Xi ;U)$, $X^U(t):=Y^U((\lambda_0^U)^{-1}(t
))$, $\tau^U$ 
defined by $(\ref{eq:extsc})$, and $\theta^U$  defined by $(\ref{eq:extscl})$, 
\[\tau^U=\lambda_0^U(\theta^U)\quad a.s..\]
\end{lemma}

\begin{proof}
It always holds 
\[\tau^U\geq\lambda_0^U(\theta^U)\quad a.s..\]
On the other hand, by Theorem \ref{th:paste}, we can suppose, 
without loss of generality, that 
\[\big(Y^U,\lambda_0^U,\Lambda_1^U,\theta^U,X^U\big)=\big(Y(\cdot
\wedge\theta ),\lambda_0(\cdot\wedge\theta ),\Lambda_1(\cdot\wedge
\theta ),\theta ,Y((\lambda_0(\cdot\wedge\theta ))^{-1}(\cdot )\wedge
\theta )\big).\]
Then, by Condition \ref{hp:ext} (iv), 
\[X^U(\tau )=Y^U((\lambda_0^U)^{-1}(\lambda_0(\theta )))=Y(\infty 
)=Y(\theta )=Y(\lambda_0^{-1}(\tau ))=X(\tau ),\quad a.s.\]
and 
\begin{eqnarray*}
X^U(\tau^{-})&=&\lim_{s\rightarrow\tau^{-}}Y((\lambda_0(\cdot\wedge
\theta ))^{-1}(s)\wedge\theta )=\lim_{s\rightarrow\tau^{-}}Y((\lambda_
0)^{-1}(s)\wedge\lambda_0^{-1}(\tau ))\\
&=&\lim_{s\rightarrow\tau^{-}}Y(\lambda_0^{-1}(s\wedge\tau ))=X(\tau^{
-})\qquad\qquad\qquad\qquad\qquad a.s.\end{eqnarray*}
Therefore 
\[\tau^U\leq\tau =\lambda_0(\theta )=\lambda_0^U(\theta^U)\quad a
.s..\]
\end{proof}
\vskip.2in

The following lemma is the analog of Theorem 3.6 of 
\cite{CK19}.

\begin{lemma}\label{th:snat}
Under Condition \ref{hp:ext}, for every solution $(Y^U,\lambda_0^
U,\Lambda_1^U)$ 
of the stopped controlled martingale problem 
for $(A,E_0,B,\Xi ;U)$, $X^U(\cdot ):=Y^U((\lambda_0^U)^{-1}(\cdot 
))$ is a natural solution of 
the stopped constrained martingale problem for 
$(A,E_0,B,\Xi ;U)$ with $\Lambda^U$ defined by $(\ref{eq:Lam})$. 
\end{lemma}

\begin{proof}
By Lemma \ref{th:extc}, $(X^U,\Lambda^U)(\cdot )=(X^U,\Lambda^U)(
\cdot\wedge\tau^U)$ 
a.s.. 

By Theorem \ref{th:paste}, we can suppose, without loss 
of generality, that 
\[\big(Y^U,\lambda_0^U,\Lambda_1^U,\theta^U,X^U\big)=\big(Y(\cdot
\wedge\theta ),\lambda_0(\cdot\wedge\theta ),\Lambda_1(\cdot\wedge
\theta ),\theta ,Y((\lambda_0(\cdot\wedge\theta ))^{-1}(\cdot )\wedge
\theta )\big).\]
Then, by Condition \ref{hp:ext} (ii), the event 
\[\{\theta^U=\infty ,\;\lim_{s\rightarrow\infty}\lambda_0^U(s)<\infty 
\}=\{\theta =\infty ,\;\lim_{s\rightarrow\infty}\lambda_0(s)<\infty 
\}\]
has zero probability. Finally, let $\{\gamma_n\}$ be the sequence of 
random variables of Condition \ref{hp:ext} (iii) and define 
\[\gamma^U_n:=\left\{\begin{array}{ll}
n&\mbox{\rm if }\theta\leq n,\\
\gamma_n&\mbox{\rm if }\theta >n.\end{array}
\right.\]
Then $\gamma^U_n\rightarrow\infty$ a.s., for each $n$ $\gamma_n^U$ is a $\bigg
\{{\cal F}^{Y^U,\lambda_0^U,\Lambda_1^U}_{(\lambda_0^U)^{-1}(t)}\bigg
\}$- stopping 
time for each $n$ and 
\[(\lambda_0^U)^{-1}(\gamma^U_n)\wedge\theta^U=\lambda_0^{-1}(\gamma^
U_n)\wedge\theta\leq n+\lambda_0^{-1}(\gamma_n)\quad a.s..\]
Therefore 
\[\begin{array}{rcl}
&&\E\bigg[\bigg|\int_{[0,t\wedge\gamma^U_n]\times\Xi}Bf(x,u)\Lambda^
U(ds\times dx\times du)\bigg|\bigg]\\
&&\leq\|Bf\|\E\big[\lambda_1^U((\lambda_0^U)^{-1}(t\wedge\gamma^U_
n)\wedge\theta^U)\big]\\
&&\leq\|Bf\|\E\big[(\lambda_0^U)^{-1}(t\wedge\gamma^U_n)\wedge\theta^
U\big]<\infty ,\\
\end{array}
\]
so that $(\ref{eq:scmp})$ is a local martingale. 
\end{proof}
\vskip.3in

\noindent {\bf Proof of Theorem \ref{th:fromst}}

First of all note that the arguments of the proof of 
point (a) of Theorem 4.2.2 of \cite{EK86} apply to 
constrained martingale problems as well, so that it is 
sufficient to prove that any two natural solutions of 
of the constrained martingale problem for $(A,E_0,B,\Xi )$ 
with the same initial distribution 
have the same one-dimensional distributions. 
The proof of Theorem 4.6.2 of \cite{EK86} essentially 
carries over. The only thing we have to check is that, 
with $V_i$ and $P_i$ as in Theorem 4.6.2 of \cite{EK86}, (in 
particular, for each $i$, $V_i=U_k{}_{}$ for some $k$) $P_i$ is the 
distribution of a natural solution of the 
stopped constrained martingale problem for $(A,E_0,B,\Xi ;V_i)$. 
To see this, let $X$ be a natural solution of 
the constrained martingale problem for $(A,E_0,B,\Xi )$. 
and suppose $X(\cdot )=Y(\lambda_0^{-1}(\cdot )$) for some solution 
$(Y,\lambda_0,\Lambda_1)$ of the 
controlled martingale problem for $(A,E_0,B,\Xi )$. (Note that 
$Y$ denotes a different object in the proof of Theorem 4.6.2 of \cite{EK86}.) 
Let 
\[\theta_0:=0,\quad\theta_i:=\inf\{t\geq\theta_{i-1}:\,Y(t)\notin 
V_i\mbox{\rm \ or }Y(t-)\notin V_i\},\,\,i\geq 1,\]
\[\rho_i:=\inf\{t\geq 0:\,Y(\theta_{i-1}+t)\notin V_i\mbox{\rm \ or }
Y((\theta_{i-1}+t)-)\notin V_i\},\mbox{\rm \ on }\{\theta_{i-1}<\infty 
\},\,\,\]
and, for $i$ such that $\P (\lambda_0(\theta_{i-1})<\infty )=\P (
\theta_{i-1}<\infty )>0$, 
\begin{eqnarray*}
&&Q_i(D):=\\
&&\frac {\E\bigg[e^{-\beta\lambda_0(\theta_{i-1})}{\bf 1}_{\{\lambda_
0(\theta_{i-1})<\infty \}}{\bf 1}_D\big(Y(\theta_{i-1}+\cdot\wedge
\rho_i),\lambda_0(\theta_{i-1}+\cdot\wedge\rho_i)-\lambda_0(\theta_{
i-1}),\Lambda_1^{\theta_{i-1},\rho_i}(\cdot )\big)\bigg]}{\E\bigg
[e^{-\beta\lambda_0(\theta_{i-1})}{\bf 1}_{\{\lambda_0(\theta_{i-
1})<\infty \}}\bigg]},\end{eqnarray*}
where $\beta$ is a positive number, 
$D\in {\cal B}(D_E[0,\infty )\times C_{[0,\infty )}[0,\infty )\times 
{\cal L}_{{\cal U}})$ and $\Lambda_1^{\theta_{i-1},\rho_i}(\cdot 
)$ is the 
measure on $[0,\infty )\times {\cal U}$ defined by 
$\Lambda_1^{\theta_{i-1},\rho_i}([0,t]\times C):=\Lambda_1([\theta_{
i-1},\theta_{i-1}+t\wedge\rho_i]\times C)$. 
Then, by Lemma 2.11 of \cite{CK19} and Remark \ref{re:tost}, 
the coordinate process on $D_E[0,\infty )\times C_{[0,\infty )}[0
,\infty )\times {\cal L}_{{\cal U}}$, 
$(\eta ,l_0,L_1)$, under $Q_i$ is a solution of the stopped controlled 
martingale problem for $(A,E_0,B,\Xi ;V_i)$ and hence, by 
Lemma \ref{th:snat}, under $Q_i$ 
$\eta ((l_0)^{-1}(\cdot ))$ is a natural solution of the stopped constrained 
martingale problem for $(A,E_0,B,\Xi ;V_i)$. 

On the other hand, with $\tau_i$ defined by $(\ref{eq:extc})$, 
and 
\[q_i:=\inf\{t\geq 0:\,X(\tau_{i-1}+t)\notin V_i\mbox{\rm \ or }X
((\tau_{i-1}+t)^{-})\notin V_i\},\]
(note that $q_i$ is denoted as $\eta^i$ in \cite{EK86}) 
we have 
\[\lambda_0^{-1}(\tau_{i-1}+t\wedge q_i)=\theta_{i-1}+\big(\lambda_
0(\theta_{i-1}+\cdot\wedge\rho_i)-\lambda_0(\theta_{i-1})\big)^{-
1}(t)\wedge\rho_i,\]
so that the distribution $P_i$ of Theorem 4.6.2 of 
\cite{EK86}, 
\[P_i(C):=\frac {\E\bigg[e^{-\beta\tau_{i-1}}{\bf 1}_{\{\tau_{i-1}
<\infty \}}{\bf 1}_C\big(X\big(\tau_{i-1}+\cdot\wedge q_i)\big)\bigg
]}{\E\bigg[e^{-\beta\tau_{i-1}}{\bf 1}_{\{\tau_{i-1}<\infty \}}\bigg
]}\]
can be written as 
\begin{eqnarray*}
&&P_i(C)=\frac {\E\bigg[e^{-\beta\lambda_0(\theta_{i-1})}{\bf 1}_{
\{\lambda_0(\theta_{i-1})<\infty \}}{\bf 1}_C\big(Y\big(\theta_{i
-1}+\big(\lambda_0(\theta_{i-1}+\cdot\wedge\rho_i)-\lambda_0(\theta_{
i-1})\big)^{-1}(\cdot )\wedge\rho_i\big)\big)\bigg]}{\E\bigg[e^{-
\beta\lambda_0(\theta_{i-1})}{\bf 1}_{\{\lambda_0(\theta_{i-1})<\infty 
\}}\bigg]}\\
&&=\E^{Q_i}\bigg[{\bf 1}_C\big(\eta\big((l_0)^{-1}(\cdot )\big)\big
)\bigg],\end{eqnarray*}
that is $P_i$ is the distribution of $\eta ((l_0)^{-1}(\cdot ))$ under $
Q_i$.  
\hfill $\Box$ \bigskip
\vskip.5in

\noindent {\bf Funding}

This research did not receive any specific grant from funding agencies in the public,
commercial, or not-for-profit sectors.
\vskip.5in

\noindent {\bf Declaration of competing interest}

The authors declare that they have no known competing financial interests or personal
relationships that could have appeared to influence the work reported in this paper.
\vskip.5in


\begin{thebibliography}{00}

\bibitem[Arendt, Kunkel and Kunze(2016)]{AKK16}
W. Arendt, S. Kunkel and M. Kunze
\newblock Diffusion with nonlocal boundary conditions. 
\newblock \emph{J. Funct. Anal.}, 270\penalty0 (7):\penalty0 2483--2507,
  2016.

\bibitem[Bass(1996)]{Bas96}
R.~F. Bass.
\newblock Uniqueness for the {S}korokhod equation with normal reflection in
  {L}ipschitz domains.
\newblock \emph{Electron. J. Probab.}, 1:\penalty0 no. 11, approx. 29 pp.\,
  1996.

\bibitem[Bass and Hsu(1991)]{BH91}
R.~F. Bass and P. Hsu.
\newblock Some potential theory for reflecting {B}rownian motion in
  {H}\"{o}lder and {L}ipschitz domains.
\newblock \emph{Ann. Probab.}, 19\penalty0 (2):\penalty0 486--508, 1991.

\bibitem[Costantini and Kurtz(2015)]{CK15}
C. Costantini and T.~G. Kurtz.
\newblock Viscosity methods giving uniqueness for martingale problems.
\newblock \emph{Electron. J. Probab.}, 20:\penalty0 no. 67, 27, 2015.

\bibitem[Costantini and Kurtz(2018)]{CK18}
C. Costantini and T.~G. Kurtz.
\newblock Existence and uniqueness of reflecting diffusions in cusps.
\newblock \emph{Electron. J. Probab.}, 23:\penalty0 no. 84, 21, 2018.

\bibitem[Costantini and Kurtz(2019)]{CK19}
C. Costantini and T.~G. Kurtz.
\newblock Markov selection for constrained martingale problems.
\newblock \emph{Electron. J. Probab.}, 24:\penalty0 no. 135, 31, 2019.

\bibitem[Costantini and Kurtz(2022)]{CK22}
C. Costantini and T.~G. Kurtz.
\newblock An inhomogeneous reverse ergodic theorem and application to a new
  uniqueness result for reflecting diffusions.
\newblock \emph{arXiv:2106.07208v4}, 2022.

\bibitem[Dai and Williams(1996)]{DW96}
J.~Dai and R.~J. Williams.
\newblock Existence and uniqueness of semimartingale reflecting brownian
  motions in convex polyhedrons.
\newblock \emph{Theory of Probability \& Its Applications}, 40\penalty0
  (1):\penalty0 1--40, 1996.

\bibitem[Davis and Norman(1990)]{DN90}
M.~H.~A. Davis and A.~R. Norman.
\newblock Portfolio selection with transaction costs. 
\newblock \emph{Math. Oper. Res.}, 15\penalty0 (4):\penalty0 676--713,
  1990.

\bibitem[DeBlassie and Toby(1993)]{DBT93}
R.~D. DeBlassie and E.~H. Toby.
\newblock Reflecting {B}rownian motion in a cusp.
\newblock \emph{Trans. Amer. Math. Soc.}, 339\penalty0 (1):\penalty0 297--321,
  1993.

\bibitem[Dupuis and Ishii(1993)]{DI93}
P. Dupuis and H. Ishii.
\newblock S{DE}s with oblique reflection on nonsmooth domains.
\newblock \emph{Ann. Probab.}, 21\penalty0 (1):\penalty0 554--580, 1993.

\bibitem[Ethier and Kurtz(1986)]{EK86}
S.~N. Ethier and T.~G. Kurtz.
\newblock \emph{Markov {P}rocesses: {C}haracterization and {C}onvergence}.
\newblock Wiley Series in Probability and Mathematical Statistics: Probability
  and Mathematical Statistics. John Wiley \& Sons Inc., New York, 1986.
\newblock ISBN 0-471-08186-8.

\bibitem[Harrison and Reiman(1981)]{HR81}
J.~M. Harrison and M.~I. Reiman.
\newblock Reflected Brownian motion on an orthant.
\newblock \emph{Ann. Probab.}, 9\penalty0 (2):\penalty0 302--308, 1981.

\bibitem[Kang and Williams(2012)]{KW12}
W.~N. Kang and R.~J. Williams.
\newblock Diffusion approximation for an input-queued switch operating under a
  maximum weight matching policy.
\newblock \emph{Stoch. Syst.}, 2\penalty0 (2):\penalty0 277--321, 2012.

\bibitem[Kang et~al.(2009)Kang, Kelly, Lee, and Williams]{KKLW09}
W.~N. Kang, F.~P. Kelly, N.~H. Lee, and R.~J. Williams.
\newblock State space collapse and diffusion approximation for a network
  operating under a fair bandwidth sharing policy.
\newblock \emph{Ann. Appl. Probab.}, 19\penalty0 (5):\penalty0 1719--1780,
  2009.

\bibitem[Kurtz(1990)]{Kur90}
T.~G. Kurtz.
\newblock Martingale problems for constrained {M}arkov problems.
\newblock In \emph{Recent advances in stochastic calculus (College Park, MD,
  1987)}, Progr. Automat. Info. Systems, pages 151--168. Springer, New York,
  1990.

\bibitem[Kurtz(1991)]{Kur91}
T.~G. Kurtz.
\newblock A control formulation for constrained {M}arkov processes.
\newblock In \emph{Mathematics of random media (Blacksburg, VA, 1989)},
  volume~27 of \emph{Lectures in Appl. Math.}, pages 139--150. Amer. Math.
  Soc., Providence, RI, 1991.

\bibitem[Kurtz and Sockbridge(2001)]{KS01}
T.~G. Kurtz and R.~H. Stockbridge.
\newblock Stationary solutions and forward equations for controlled and singular martingale problems. 
\newblock \emph{Electron. J. Probab.}, 6 :\penalty0 no. 17, 52, 2001.

\bibitem[Kwon and Williams(1991)]{KW91}
Y.~Kwon and R.~J. Williams.
\newblock Reflected {B}rownian motion in a cone with radially homogeneous
  reflection field.
\newblock \emph{Trans. Amer. Math. Soc.}, 327\penalty0 (2):\penalty0 739--780,
  1991.

\bibitem[Menaldi and Robin(1985)]{MR85}
J-L. Menaldi and M. Robin.
\newblock Reflected diffusion processes with jumps. 
\newblock \emph{Ann. Probab.}, 13\penalty0 (2):\penalty0 319--341,
  1985.

\bibitem[Reiman and Williams(1988)]{RW88}
M.~I. Reiman and R.~J. Williams.
\newblock A boundary property of semimartingale reflecting {B}rownian motions.
\newblock \emph{Probab. Theory Related Fields}, 77\penalty0 (1):\penalty0
  87--97, 1988; 80\penalty0 (4):\penalty0 633, 1989.

\bibitem[Saisho(1987)]{Sai87}
Y. Saisho.
\newblock Stochastic differential equations for multidimensional domain with
  reflecting boundary.
\newblock \emph{Probab. Theory Related Fields}, 74\penalty0 (3):\penalty0
  455--477, 1987.

\bibitem[Shreve and Soner(1994)]{SS94}
S.~E. Shreve and H.~M. Soner.
\newblock Optimal investment and consumption with transaction costs. 
\newblock \emph{Ann. Appl. Probab.}, 4\penalty0 (3):\penalty0 609--692,
  1994.

\bibitem[Tanaka(1979)]{Tan79}
H. Tanaka.
\newblock Stochastic differential equations with reflecting boundary condition
  in convex regions.
\newblock \emph{Hiroshima Math. J.}, 9\penalty0 (1):\penalty0 163--177, 1979.

\bibitem[Taylor and Williams(1993)]{TW93}
L.~M. Taylor and R.~J. Williams.
\newblock Existence and uniqueness of semimartingale reflecting {B}rownian
  motions in an orthant.
\newblock \emph{Probab. Theory Related Fields}, 96\penalty0 (3):\penalty0
  283--317, 1993.

\bibitem[Varadhan and Williams(1985)]{VW85}
S.~R.~S. Varadhan and R.~J. Williams.
\newblock Brownian motion in a wedge with oblique reflection.
\newblock \emph{Comm. Pure Appl. Math.}, 38\penalty0 (4):\penalty0 405--443,
  1985.

\bibitem[Williams(1985)]{Wil85}
R.~J. Williams.
\newblock Reflected {B}rownian motion in a wedge: semimartingale property.
\newblock \emph{Z. Wahrsch. Verw. Gebiete}, 69\penalty0 (2):\penalty0 161--176,
  1985.



\end{thebibliography}

\end{document}